\documentclass[12pt,a4paper]{article}

\usepackage{amsmath,amssymb,amsthm,fancyhdr,supertabular,calc,eepic,fullpage}
\usepackage{mathrsfs}
\usepackage[mathcal]{euscript}
\usepackage{wasysym}

\usepackage[dvips]{graphicx}
\usepackage[dvips]{color}

\newtheorem{theorem}{Theorem}
\newtheorem{lemma}[theorem]{Lemma}

\newtheorem{corollary}[theorem]{Corollary}

\newtheorem*{theorema}{Theorem A}

\newtheorem*{hypothesis}{Hypothesis}
\newtheorem*{remark}{Remark}

\newcommand{\GL}{G_\linel}
\newcommand{\ordGL}{|G_\linel|}
\newcommand{\Galph}{G_\alpha}

\newcommand{\ordGalph}{|G_\alpha|}

\newcommand{\point}{\Pi}
\newcommand{\lines}{\Lambda}
\newcommand{\linel}{\mathfrak{L}}
\newcommand{\hatt}{\hat {\ } }
\newcommand{\spaceS}{\mathcal{S}}

\begin{document}


\pagestyle{plain}

\pagenumbering{arabic}
\title{$PSL(3,q)$ and line-transitive linear spaces}
\author{Nick Gill\thanks{This paper contains results from the
    author's PhD thesis; I would like to thank my supervisor Professor Jan Saxl. I would also like to thank Professor Francis Buekenhout for his excellent advice.}}

\maketitle

\begin{abstract}
We present a partial classification of those finite linear spaces $\spaceS$ on
which an almost simple group $G$ with socle $PSL(3,q)$ acts line-transitively.
\end{abstract}

A {\it linear space} $\spaceS$ is an incidence structure consisting of
a set of points $\point$ and a set of lines $\lines$ in the power set
of $\point$ such that any two points are incident with exactly one
line. The linear space is called {\it non-trivial} if every line
contains at least three points and there are at least two lines. Write
$v=|\point|$ and $b=|\lines|$.

The investigation of those finite linear spaces which admit an almost simple group 
that is transitive upon lines is already underway \cite{cnp, camspiez}
motivated largely by the theorem of Camina and Praeger \cite{campraeger}. We
continue this investigation by considering the situation
when the socle of a line-transitive automorphism group is $PSL(3,q)$. The
statement of our theorem is as follows:

\begin{theorema}
Suppose that $PSL(3,q)\unlhd G \leq Aut PSL(3,q)$ and that $G$ acts
line-transitively on a finite linear space $\spaceS.$ Then one of the
following holds:
\begin{itemize}
\item $S=PG(2,q)$, the Desarguesian projective plane, and $G$ acts
  2-transitively on points;
\item $PSL(3,q)$ is point-transitive but not line-transitive on
  $\spaceS$. Furthermore, if $\Galph$ is a point-stabilizer in $G,$ then
  $\Galph\cap PSL(3,q)\cong PSL(3,q_0)$ where $q=q_0^a$ for some integer $a$.
\end{itemize}
\end{theorema}

The second possibility given in Theorem A is not expected to yield any examples, however it has not yet been excluded. The proof of Theorem A will depend heavily upon a result of Camina, Neumann and Praeger which classifies the line-transitive actions of $PSL(2,q)$:

\begin{theorem}\label{theorem:psl2}
Let $G = PSL(2,q), q\geq 4$ and suppose that $G$ acts
line-transitively on a linear space $\spaceS$.  Then one of the
following holds: 
\begin{itemize}
\item $G=PSL(2,2^a), a\geq 3$ acting transitively on $\spaceS$, a
Witt-Bose-Shrikhande space. Here $\point$ is the set of dihedral
subgroups of $G$ of order $2(q+1)$ and $\lines$ is the set of involutions
$t\in G$ with the incidence relation being inclusion.

\item $\spaceS=PG(2,2)$, $G=PSL(2,7)$ and the action is 2-transitive.
\end{itemize}
\end{theorem}

Theorem \ref{theorem:psl2} has not appeared in the literature. A weaker version of the result has been proven by Liu \cite{liupsl2}; at the end of Section \ref{subsection:ltls} we will prove Theorem \ref{theorem:psl2} using Liu's result.

In the case where $\spaceS$ is a projective plane Theorem A is implied by \cite[Theorem A]{gill}. Furthermore where $G$ is flag-transitive upon $\spaceS$ Theorem A is implied by \cite{saxl}.

Now observe that if a linear space $\spaceS$ is line-transitive then
every line has the same number, $k$, of points and every point lies on the same
number, $r$, of lines. Such a linear space is called {\it regular} and 
those line-transitive linear spaces for which $k=3$ or $k=4$ have been
completely classified in \cite{clapham, ks, ftclass,
camsiem, li}. 

Hence in order to prove Theorem A we need to consider the situation when $\spaceS$
is not a projective plane, is not flag-transitive and $k\geq 5$. The rest of the paper will be occupied
with this proof. The first two sections outline some background
lemmas concerning linear spaces. Section \ref{section:background}
gives background information about $PSL(3,q)$. In Section
\ref{section:simplereduction} we introduce the notion of {\bf
  exceptionality} of permutation representations, the relevance of which was
pointed out by Professor Peter Neumann. In Section \ref{section:proof} we start the proof proper; we set up a hypothesis under which the proof can proceed, we cover some preliminary cases and we split the remainder of the proof into several cases. In the remaining sections we examine these cases one at a time.

The following notation will hold, unless stated otherwise, throughout
this paper. We will take $G$ to be a group acting on a regular
linear space $\spaceS$ with parameters $b,v,k,r$. We will write
$\alpha$ to be a point of $\spaceS$ with $G_\alpha$ to be the
stabilizer of $\alpha$ in the action of $G$. Similarly $\linel$ is a
line of $\spaceS$ and $\GL$ is the corresponding line-stabilizer.

\section{Known Lemmas}\label{section:knownlemmas}

We list here some well-known lemmas which we will use later. The first lemma is proved easily by counting.

\begin{lemma}\label{lemma:arith}
\begin{enumerate}
\item	$b=\frac{v(v-1)}{k(k-1)}\geq v$ (Fisher's inequality);
\item	$r=\frac{v-1}{k-1}\geq k$;
\end{enumerate}
\end{lemma}

\begin{lemma}\label{lemma:bdoublev}\cite[Lemma 6.5]{cnp}
Let $p$ be an odd prime divisor of $v$. 
\begin{enumerate}
\item	If $b=\frac{3}{2}v$ then $p=5$ and $25\not\big| v$, or $p\equiv 1, 2, 4$ or $8(15)$;
\item	If $b=2v$ then $p\equiv 1(4)$.
\end{enumerate}
\end{lemma}

For the remainder of this section assume that $G$ acts
line-transitively on the linear space $\spaceS$.

\begin{theorem}\label{theorem:kdivv}\cite[Theorem1]{camgag}
If $k\big|v$ then $G$ is flag-transitive.
\end{theorem}

\begin{lemma}\label{lemma:invft}\cite[Lemma 4]{camsiem}
If $g$ is an involution of $G$ and $g$ fixes no points, then $k\big|v$. In particular, $G$ is flag-transitive.
\end{lemma}

\begin{lemma}\label{lemma:fixedset}\cite[Lemma 2]{camsiem}
Let $\linel$ be a line in $\spaceS$ and let $T\leq \GL$. Assume that $T$ satisfies the following two conditions:
\begin{enumerate}
\item	$|Fix_\point(T)\cap \linel|>1$;
\item	if $U\leq \GL$ and $|Fix_\point(U)\cap L|>1$ and $U$ is conjugate to $T$ in $G$, then $U$ is conjugate to $T$ in $\GL.$
\end{enumerate}
Then either $Fix_\point(T)\subseteq \linel$ or the induced linear space on $Fix_\point(T)$ is regular and $N_G(T)$ acts line-transitively on the space.
\end{lemma}

\begin{lemma}\label{lemma:indextwo}\cite[Lemma 2.2]{camspiez}
Let $g$ be an involution in $G$ and assume that there exists $N$, $N\lhd G$ such that $|G:N|=2$ with $g\not\in N$. Then $N$ acts line-transitively also.
\end{lemma}

Note that Lemma \ref{lemma:indextwo} allows us to conclude that if
$PGL(2,q)$ acts transitively on the lines of a linear space $\spaceS$ then
$PSL(2,q)$ also acts transitively on the lines of $\spaceS$ and so that space is known.

Our next result provides the framework for our analysis of the
line-transitive actions
of $PSL(3,q)$. Since $\spaceS$ is not a projective plane then, by
Fisher's inequality $b>v$ and since $b=v(v-1)/(k(k-1))$, there must be
some prime $p$ that divides both $v-1$ and $b$. We shall refer to such
a prime as a {\it significant} prime. 

\begin{lemma}\label{lemma:sigprime}\cite[Lemma 6.1]{cnp}
Suppose that $\spaceS$ is not a projective plane and let $p$ be a
significant prime. Let $P$ be a Sylow $p$-subgroup of $G_\alpha$. Then $P$ is a Sylow $p$-subgroup of $G$ and $G_\alpha$ contains the normalizer $N_G(P)$.
\end{lemma}

\begin{lemma}\label{lemma:cong}\cite[Lemma 6.3]{cnp}
Let $H,K$ be subgroups such that 
$$G_\alpha\leq H<K\leq G$$
and let $c=|K:H|$. Then $r$ divides $\frac{1}{2}(c-1)k$ and $b$ divides $\frac{1}{2}(c-1)v$.
\end{lemma}

\begin{corollary}\label{corollary:supergroups}\cite[Corollary 6.4]{cnp}
Let $H, K$ be as in Lemma \ref{lemma:cong}.
\begin{enumerate}
\item	Let
	$$c_0= {\rm gcd}\{(c-1) \ | \ c=|K:H|, \ where \ G_\alpha\leq H<K\leq G\}.$$
	Then $r$ divides $\frac{1}{2}c_0k$ and $b$ divides $\frac{1}{2}c_0v$.
\item	There cannot be groups $H,K$ such that $G_\alpha\leq H<K\leq G$ and $|K:H|=2$.
\item	If there are groups $H,K$ such that $G_\alpha\leq H<K\leq G$ and $|K:H|=3$ then $\spaceS$ is a projective plane.
\end{enumerate}
\end{corollary}

\section{New Lemmas}\label{section:newlemmas}

We state a series of lemmas which will be used in our analysis of the
actions of $PSL(3,q)$. The first is a generalization of the Fisher
inequality to non-regular linear spaces.

\subsection{General linear spaces}

\begin{lemma}\label{lemma:generalfisher}
In any linear space $\spaceS$, not necessarily regular, Fisher's inequality holds: $b\geq v$.
\end{lemma}
\begin{proof}
We need to prove the statement under the assumption that the number of points in a line is not a constant. Let $c$ be the maximum number of points on a line of $\spaceS$. Since any two points lie on a unique line we know that
$$b\geq\frac{{v\choose2}}{{c\choose2}}=\frac{v(v-1)}{c(c-1)}.$$
Thus if $c(c-1)\leq v-1$ then we are finished. Assume to the contrary from this point on. We split into two cases:

\begin{enumerate}
\item	{\bf Suppose that $(c-1)^2\geq v$.} Then $\frac{v+c-1}{c}\leq
c-1$. Let $\linel_c$ be a line with $c$ points on it and choose $\alpha$ a
point not on $\linel_c$. Then the average number of points in a line
containing $\alpha$ and intersecting $\linel_c$ is less than or equal to
$\frac{v-c-1}{c}+2=\frac{v+c-1}{c}\leq c-1.$ Call the number of lines
intersecting $\linel_c$, $b_0$, and observe that
$$b\geq b_0\geq\frac{(v-c)c}{c-1}.$$
Now we know that $v>c^2$ and so $vc-c^2>vc-v$, hence
$\frac{(v-c)c}{c-1}>v$. Thus this case is covered.
\item	{\bf Suppose that $(c-1)^2<v\leq c(c-1)$.} Note that $v>2$
  implies that $c>2$. Let $r_\alpha$ be the number of lines incident with a point $\alpha$. If $r_\alpha\geq c$ for all $\alpha$ then, let $f$ be the number of flags:
$$vc\leq f\leq bc.$$
Thus $v\leq b$ as required. Assume then that there exists a point
$\alpha $ such that $r_\alpha\leq c-1$. Observe that every line not
passing through $\alpha$ must have be incident with at most $r_\alpha$
points. Remove $\alpha$ and any lines which are incident with only
$\alpha$ and one other point. Then $v>c$ and we
still have a linear space, $S^*$. $S^*$ has $v-1$ points, at most $b$
lines, and the maximum number of points on a line is $c-1$. This
implies that,
$$b_S\geq b_{S^*}\geq \frac{{{v-1}\choose2}}{{{c-1}\choose2}}=\frac{(v-1)(v-2)}{(c-1)(c-2)}.$$
Thus we are finished so long as $(c-1)(c-2)< v-2$. But $(c-1)^2<v$
gives us this inequality since $c \geq 3$.
\end{enumerate}

All cases are proved and the result stands.
\end{proof}

\subsection{Regular linear spaces}

We return to our assumption that $\spaceS$ is a regular linear space, recall that this means that all lines are incident with the same number of points.

\begin{lemma}\label{lemma:invfix}
Let $g\in G$ be an involution. Then $g$ fixes at least $(v-1)/k$ lines. 
\end{lemma}
\begin{proof}
If $g$ has no fixed point then $g$ fixes $v/k\geq (v-1)/k$ lines. If
$g$ has a fixed point, $\alpha$, then let $m$ be the number of fixed
lines through $\alpha$. By definition, $g$ moves the rest of the lines
through $\alpha$. Apart from $\alpha$ these lines contain $v-m(k-1)-1$
points. None of these points is fixed hence every one of these points
lies on a fixed line. Thus the number of lines fixed by $g$ is at
least
$$m+\frac{v-m(k-1)-1}{k}=\frac{v+m-1}{k}\geq \frac{v-1}{k}$$
lines as required.
\end{proof}

\begin{lemma}\label{lemma:fixpointline}
Let $g$ be an involution which is an automorphism of a linear space
$\spaceS$. Suppose that $\spaceS$ has a constant number of points on a
line, $k$, and that $g$ fixes $d_l$ lines and $d_p$ points. Then,
either
\begin{itemize}
\item	$d_l\geq d_p$; or
\item	$v=k^2$.
\end{itemize}
\end{lemma}
\begin{proof}
We know that if $\spaceS$ is a projective plane then the result holds
since the permutation character on points and lines is the same
\cite[4.1.2]{dembov}. Now suppose that $\spaceS$ is not a projective
plane and split into two cases:
\begin{enumerate}
\item	{\bf Suppose that $d_p\leq k$.} Assume that $d_l<d_p$. We know, by Lemma \ref{lemma:invfix}, that $g$ fixes at least $\frac{v-1}{k}$ lines. Then
\begin{eqnarray*}
d_l<d_p&\implies&\frac{v-1}{k}<k \\
&\implies& v-1<k^2. 
\end{eqnarray*}
Then, since $(k-1)\big|(v-1)$, we must have $\frac{v-1}{k-1}\leq
k+1$. If $\frac{v-1}{k-1}\leq k$ then $b\leq v$ and so $b=v$ and $\spaceS$ is a projective plane. If $\frac{v-1}{k-1}= k+1$ then $v=k^2$ as given.
\item	{\bf Suppose that $d_p>k$.} Then the fixed points and lines of
  $g$ form a linear space. We may appeal to Lemma \ref{lemma:generalfisher}.
\end{enumerate}
\end{proof}



\begin{lemma}\label{lemma:sigprimesmall}
Suppose that $b=\frac{c}{d}v$ where $(c,d)=1$. Then the significant primes are exactly those which divide $c$.
\end{lemma}
\begin{proof}
By definition a prime is significant if it divides $b$ and $v-1$. Then we just use the fact that
$$\frac{c}{d}v=b=\frac{v(v-1)}{k(k-1)}=\frac{(v-1)/(k-1)}{k}v.$$
\end{proof}

\begin{lemma}\label{lemma:normalizersub}
Let $H<\Galph$. If $N_G(H)\not\leq \Galph$ then $H$ is in $\GL$ for some line $\linel$.
\end{lemma}
\begin{proof}
Simply take $g\in N_G(H)\backslash \Galph$. Then $H^g=H$ is contained in $\Galph$ and $G_{\alpha g}$. Hence $H$ fixes the line joining $\alpha$ and $\alpha g$.
\end{proof}

\subsection{Line-transitive linear spaces}\label{subsection:ltls}

Throughout this section we assume that $G$ acts line-transitively on $\spaceS$.

\begin{lemma}\label{lemma:invineq}
Let $g$ be an involution of $G$ and write $n_g=|g^G|$ for the size of a
conjugacy class of involutions in $G$. Let $r_g=|g^G\cap \GL|$ be the
number of such involutions in a line-stabilizer $\GL$. Then the
following inequality holds:
$$\frac{n_g(v-1)}{br_g}\leq k\leq\frac{r_gv}{n_g}+1.$$
\end{lemma}
\begin{proof}
Count pairs of the form $(\linel,g)$ where $\linel$ is a line and $g$ is an involution fixing $\linel$, in two different ways. Then
$$|\{(\linel,g)\}|=br_g\geq n_gc$$
where $c$ is the minimum number of lines fixed by an involution. Now, by the previous lemma, $c\geq \frac{v-1}{k}$ thus we have
$$r_g\geq\frac{n_gc}{b}\geq \frac{n_g(v-1)}{bk}=\frac{n_g(k-1)}{v}.$$
This implies two inequalities:
$$k-1\leq\frac{r_gv}{n_g}\quad,\quad k\geq\frac{n_g(v-1)}{br_g}$$
and the result follows.
\end{proof}

\begin{lemma}
Suppose that $|G_\alpha|=\frac{c}{d}|\GL|$ where $(c,d)=1$. Then the significant primes are exactly those which divide $c$.
\end{lemma}
\begin{proof}
Simply use the fact that $v=|G|/|G_\alpha|$, $b=|G|/|\GL|$ and refer
to Lemma \ref{lemma:sigprimesmall}.
\end{proof}

%

\begin{lemma}\label{lemma:largeprime}
Suppose that $p^a$ is a prime power dividing $v-1$ and that $p$ does not divide into $|G|$. Then $p^a$ divides $k(k-1)$.
\end{lemma}
\begin{proof}
Since $p$ does not divide $|G|$, $p$ cannot divide into $b$. Since $b=\frac{v(v-1)}{k(k-1)}$ and $p^a$ divides into $v-1$ we must have $p^a$ dividing into $k(k-1)$.
\end{proof}

We will often repeatedly use Lemma \ref{lemma:largeprime}, with
different primes, to exclude the
possibility of a particular group, $G,$ acting line-transitively on a
space with a particular number of points, $v$. Our method for doing
this usually involves showing that any line size $k$ must be too large
to satisfy Fisher's inequality (Lemma \ref{lemma:arith}). 

As promised we conclude this section with a proof of Theorem \ref{theorem:psl2}:

\begin{proof}
We suppose that $G=PSL(2,q)$ acts line-transitively on a linear space $\spaceS$ and seek to demonstrate that only the possibilities listed in Theorem \ref{theorem:psl2} can occur. We take as our starting point the main result of \cite{liupsl2} which gives two extra possibilites which we must exclude.

Firstly if $\spaceS$ is a projective plane then \cite[Theorem A]{gill} implies that $G=PSL(2,7)$ and $\spaceS=PG(2,2)$ as required.

The other possibility which Liu leaves open is that $q=2^6$ and $(v,k)=(2080, 12)$. In this case a line-stabilizer in $G$ has size $8$ and so contains at most 7 involutions. Applying Lemma \ref{lemma:invineq} we have
$$k\leq\frac{r_gv}{n_g}+1\leq\frac{2080\times 7}{63\times 65}+1<5.$$
This is a contradiction and this possibility can be excluded as required.
\end{proof}

\section{Background Information on $PSL(3,q)$}\label{section:background}

We will sometimes precede the structure
of a subgroup of a projective group with $\hatt$ which means that we
are giving the structure of the pre-image in the corresponding linear
group. We will also refer to elements of this linear group in terms of
matrices under the standard modular representation. 

\subsection{Subgroup information}

We need information about the subgroups of $PSL(3,q)$, $PSL(2,q)$ and
$GL(2,q)$. Here $q=p^a$ for a prime $p$, positive integer $a$.

\begin{theorem}\cite{kleidman, mitchell, bloom, hartley}\label{theorem:psl3subgroups}
The maximal subgroups of $PSL(3,q)$ are among the following
list. Conditions given are necessary for existence and maximality but
not sufficient. The first three types are all maximal for $q\geq 5$.

\begin{center}
\begin{tabular}{c|c|c}
&{\bf Description} & {\bf Notes} \\
\hline
1&$\hatt[q^2]:GL(2,q)$ & two $PSL(3,q)$-conjugacy classes \\
2&$\hatt(q-1)^2:S_3$ & one $PSL(3,q)$-conjugacy class \\
3&$\hatt (q^2+q+1).3$ & one $PSL(3,q)$-conjugacy class \\
4&$PSL(3,q_0).(q-1,3,b)$ & $q=q_0^b$ where $b$ is prime \\
5&$PSU(3,q_0)$ & $q=q_0^2$ \\
6&$A_6$ &\\
7&$3^2.SL(2,3)$ & $q$ odd \\
8&$3^2.Q_8$ & $q$ odd \\
9&$SO(3,q)\cong PGL(2,q)$ & $q$ odd \\
10&$PSL(2,7)$ & $q$ odd
\end{tabular}
\end{center}
\end{theorem}

We will refer to maximal subgroups of $PSL(3,q)$ as being of type
$x$, where $x$ is a number between 1 and 10 corresponding to the list above. 

Referring to \cite{bloom, kleidman} we state the following lemma:
\begin{lemma}\label{lemma:psl3subgroups}
Suppose that $H$ is a subgroup of $PSL(3,q)$ lying in a maximal
subgroup of type 4 or 5 and $H$ does not lie in any other maximal
subgroup of $PSL(3,q)$. Then one of the following holds:
\begin{itemize}
\item $H$ has a cyclic normal subgroup of index less than or equal to
3;
\item $H$ has socle $PSL(3,q_1), PSU(3,q_1)$ or $PSL(2,q_1)$  with index less than or equal to 3. Here $q=q_1^c$, $c$ an integer;
\item $H$ is a subgroup of $3^2.Q_8$ or $3^2.SL(2,3)$;
\item $H$ is isomorphic to $A_5$ or $PSL(2,7)$;
\item $H$ is isomorphic to $A_6$ and $q\equiv 1$ or $19(30)$;
\item $H$ is isomorphic to $A_6$, $A_6.2$ or $A_7$, furthermore $q=5^a, a$ even.
\end{itemize}
\end{lemma}

We state a result given by Suzuki \cite[Theorem 6.25]{suzuki} which gives the structure of all the subgroups of $PSL(2,q)$: 

\begin{theorem}\label{theorem:pslsubgroup}
Let $q$ be a power of the prime $p$. Let $d=(q-1,2)$. Then a subgroup of $PSL(2,q)$ is isomorphic to one of the following groups.
\begin{enumerate}
\item	The dihedral groups of order $2(q\pm 1)/d$ and their subgroups.
\item	A parabolic group $P_1$ of order $q(q-1)/d$ and its subgroups. A Sylow $p$-subgroup $P$ of $P_1$ is elementary abelian, $P\lhd P_1$ and the factor group $P_1/P$ is a cyclic group of order $(q-1)/d$.
\item	$PSL(2,r)$ or $PGL(2,r)$, where $r$ is a power of $p$ such that $r^m =q$.
\item	$A_4, S_4$ or $A_5$.
\end{enumerate}
\end{theorem}

Note that when $p=2$, the above list is complete without the final
entry. Furthermore, referring to \cite{kleidman}, we see that there
are unique $PSL(2,q)$ conjugacy classes of the maximal dihedral
subgroups of size $2(q\pm 1)/d$ as well as a unique $PSL(2,q)$
conjugacy class of parabolic subgroups $P_1$.

We will also need the subgroups of $GL(2,q)$ which can be easily obtained from
the subgroups of $PSL(2,q)$ (for the odd characteristic case see \cite[Theorem
3.4]{bloom}.)

\begin{theorem}\label{theorem:gl2subgroups}
$H$, a subgroup of $GL(2,q)$, $q=p^a$, is amongst the following up to
conjugacy in $GL(2,q)$. Note that the last two cases may be omitted
when $p=2$.
\begin{enumerate}
\item	$H$ is cyclic;
\item	$H =AD$ where
	$$A\leq\left\{\left(\begin{array}{cc} 1 & 0 \\ \lambda & 1 \\ \end{array}\right):\lambda\in GF(q)\right\}$$
	and $D\leq N(A)$, is a subgroup of the group of diagonal matrices;
\item	$H=\langle c, S\rangle$ where $c|q^2-1$, $S^2$ is a scalar 2-element in $c$;
\item	$H=\langle D,S\rangle$ where $D$ is a subgroup of the group of diagonal matrices, $S$ is an anti-diagonal 2-element and $|H:D|=2$;
\item	$H=\langle SL(2,q_0),V\rangle$ or contains $\langle SL(2,q_0),V\rangle$ as a subgroup of index 2
and here $q=q_0^c$, $V$ is a scalar matrix. In the second case, $q_0>3$;
\item	$H/\langle-I\rangle$ is isomorphic to $S_4\times C$, $A_4\times C$, or (with $p\neq 5$) $A_5\times C$, where $C$ is a scalar subgroup of $GL(2,q)//\langle-I\rangle$; 
\item	$H//\langle-I\rangle$ contains $A_4\times C$ as a subgroup of index 2 and $A_4$ as a subgroup with cyclic quotient group, $C$ is a scalar subgroup of $GL(2,q)//\langle-I\rangle$.
\end{enumerate}
\end{theorem}

We will refer to maximal subgroups of $GL(2,q)$ as being of type
$x$, where $x$ is a number between 1 and 7 corresponding to the list above.

Finally observe that $PSL(3,q)$ contains a single conjugacy class
of involutions. This class is of size $q^2(q^2+q+1)$ for $q$ odd and
of size $(q^2-1)(q^2+q+1)$ for $q$ even. In addition note that we will write
$\mu$ for $(q-1,3)$.  

\subsection{The subgroup $D$}

We define $D$ to be the centre of a Levi complement of a particular
parabolic subgroup. Typically $D$ is the projective image of  
$$\left\{\left(\begin{array}{ccc} \frac{1}{a^2} & 0
& 0 \\ 0 & a & 0 \\ 0 & 0 & a \\ \end{array}\right):a\in \mathbb{F}_q^*\right\}.$$

Suppose that $G=PSL(3,q)$ acts line-transitively on a linear
space. Since $D$ normalizes a Sylow $t$-subgroup of $PSL(3,q)$ for
many different $t$, $D$ often lies inside a point-stabilizer
$\Galph$. Furthermore, since $D$ has a large normalizer, $\hatt
GL(2,q)$, by Lemma \ref{lemma:normalizersub}, $D$ often lies inside a
line-stabilizer, $\GL$. 

We exploit this fact using Lemma
\ref{lemma:fixedset} since if $D$ satisfies the conditions given in
the lemma and the fixed points of $D$ are not collinear then we induce
a line-transitive action of $PGL(2,q)$ on a linear space. All such
actions on a non-trivial linear space are 
known. In the event that the fixed set is a trivial linear space (that
is, $k=2$) line-transitivity is equivalent to $2$-homogeneity on
points and these actions are also all well-known.

We need information about the occurrence of $D$ in various subgroups
and about how $G$-conjugates of $D$ intersect. We state the relevant
facts below; proofs are omitted as the results are easily derived from matrix
calculations.

\begin{lemma}\label{lemma:dconjintersect}
The $PSL(3,q)$ conjugates of $D$ intersect trivially.
\end{lemma}

\begin{lemma}\label{lemma:dconjugates}
Let $U:\hatt GL(2,q)$ be  a parabolic subgroup of $PSL(3,q)$, $q>7$, $U$ an elementary abelian $p$-group. We can choose $\hatt GL(2,q)$ conjugate to
$$C_G(D)=\hatt\left\{\left(\begin{array}{ccc} \frac{1}{DET} & 0 & 0\\ 0 & e & f \\
0 & g & h \\ \end{array}\right): \left(\begin{array}{cc} e & f \\ g &
h \\ \end{array}\right)\in GL(2,q), DET=eh-fg \right\}.$$

Let $H$ be a maximal subgroup of $\hatt GL(2,q)$ in $PSL(3,q)$. Write $a$ for a
primitive element of $GF(q)$.

\begin{enumerate}
\item	If $H$ is of type 2 in $\hatt GL(2,q)$ then some $\hatt
GL(2,q)$ conjugate of $H$ contains one individual conjugate, and two
families of conjugates, of $D$, generated by the projective images of
the following matrices, for $f\in GF(q)$:
$$\left(\begin{array}{ccc} \frac{1}{a^2} & 0 & 0\\ 0 & a & 0 \\ 0 & 0 & a \\ \end{array}\right), \ 
\left(\begin{array}{ccc} a & 0 & 0\\ 0 & \frac{1}{a^2} & 0 \\ 0 & f & a \\ \end{array}\right), \
\left(\begin{array}{ccc} a & 0 & 0\\ 0 & a & 0 \\ 0 & f & \frac{1}{a^2} \\ \end{array}\right).$$

\item	If $H$ is of type 3 in $\hatt GL(2,q)$ then $H$ contains only $D$.

\item	If $H$ is of type 4 in $\hatt GL(2,q)$ then some $\hatt
GL(2,q)$-conjugate of $H$ contains three conjugates of $D$, generated by the projective images of the following matrices:
$$\left(\begin{array}{ccc} \frac{1}{a^2} & 0 & 0\\ 0 & a & 0 \\ 0 & 0 & a \\ \end{array}\right), \ 
\left(\begin{array}{ccc} a & 0 & 0\\ 0 & \frac{1}{a^2} & 0 \\ 0 & 0 & a \\ \end{array}\right), \
\left(\begin{array}{ccc} a & 0 & 0\\ 0 & a & 0 \\ 0 & 0 & \frac{1}{a^2} \\ \end{array}\right).$$

\item	If $H$ is of type 5 in $\hatt GL(2,q)$ then one of the
following holds:
	\begin{itemize}
	\item $H$ contains only $D$; 
	\item $H\geq SL(2,q)$; 
	\item $H\geq SL(2,q_0)$ where $q=q_0^2$ and $q_0=3,4$ or $7$.
\end{itemize}

\item	If $H$ is of type 6 or 7 in $\hatt GL(2,q)$ then one of the
following holds:
	\begin{itemize}
	\item	$H$ contains only the central copy of $D$;
	\item	$q=13, 16$ or $19$.
	\end{itemize}

\end{enumerate}
\end{lemma}

\begin{corollary}
A subgroup of $PSL(3,q)$ of type 3 contains only the 3 diagonal
conjugates of $D$ as listed above for $H$ of type 4 in $\hatt
GL(2,q)$.
\end{corollary}

\section{Exceptional analysis}\label{section:simplereduction}

Let $B$ be a normal subgroup in a group $G$ which acts upon a set $\Pi$. Then
$(G,B,\Pi)$ is called {\it exceptional} if the only
common orbital of $B$ and $G$ in their action upon $\Pi$ is the diagonal (see \cite{gms}). 

\begin{lemma}\label{lemma:exceptional}
Suppose a group $G$ acts line-transitively on a linear space
$\spaceS$; suppose furthermore that $B$ is a normal subgroup in $G$
which is not line-transitive on $\spaceS$; finally suppose that
$|G:B|=t$, a prime. 

Then either $\spaceS$ is a projective plane or $(G,B,\Pi)$ is exceptional.
\end{lemma}
\begin{proof}
The suppositions mean that, for a line $\linel$ of $\spaceS$,
$\GL=B_\linel$. We have two possibilities:

\begin{itemize}
\item  {\bf Suppose that $B$ is point-transitive on $\spaceS$.} Then let
  $\alpha$ and $\beta$ be members of $\point$, the set of points of
  $\spaceS$. Let $\linel$ be the line connecting them. Then, since
  $G_{\alpha,\beta}\leq \GL$ and $B_{\alpha,\beta}\leq B_\linel$, we
  know that $G_{\alpha,\beta}=B_{\alpha,\beta}$.

We know furthermore that $|\Galph:B_\alpha|=t$, hence we may conclude
that, for all pairs of points $\alpha$ and $\beta$,
$|B_\alpha:B_{\alpha,\beta}|<|\Galph:G_{\alpha,\beta}|$. In other words
$(G,B,\Pi)$ is {\bf exceptional}. 

\item  {\bf Suppose that $B$ is not point-transitive on $\spaceS$.}
  Then, by the Frattini argument, $G=N_G(P)B$ for all $P\in Syl_pB$
  where $p$ is any prime dividing into $|B|$. If $\Galph\geq N_G(P)$
  then $B$ is point-transitive which is a contradiction. Thus, by
  Lemma \ref{lemma:normalizersub}, if a Sylow $p$-subgroup of $B$
  stabilizes a point then it also stabilizes a line.

Now let $b_B=|B:B_\linel|$, $v_B=|B:B_\alpha|$. Then primes dividing into
$b_B$ are a subset of the primes dividing into $v_B$. Furthermore $b=tb_B$
and $v=tv_B$. Thus primes dividing into $b$ are a subset of the primes
dividing into $v$. Thus there are no significant primes and $\spaceS$
is a projective plane.
\end{itemize}
\end{proof}

\section{The Proof of Theorem A}\label{section:proof}

In this section we commence the proof of Theorem A. We begin, in Section \ref{subsection:simplicity}, by reducing the proof to the case where we have $PSL(3,q)$ acting line-transitively on a linear space $\spaceS$. In Sections \ref{subsection:easy1} and \ref{subsection:easy2} we rule out some easy preliminary cases. Finally, in Section \ref{subsection:table}, we split the remainder of the proof into different cases to be examined one at a time in the succeeding sections.

\subsection{Reduction to Simplicity}\label{subsection:simplicity}

\begin{lemma}
Suppose that $PSL(3,q)\unlhd G\leq Aut PSL(3,q)$ and $G$ acts line-transitively on a finite linear space $\spaceS$. One of the following must hold:
\begin{itemize}
\item 	$PSL(3,q)$ acts line-transitively on $\spaceS$;
\item		The second possibility given in Theorem A holds.
\end{itemize}
\end{lemma}
\begin{proof}
Recall that, by \cite[Theorem A]{gill}, we may assume that $\spaceS$ is not a projective plane. Suppose that $PSL(3,q)$ is not line-transitive on $\spaceS$. Then there exist groups $G_1,G_2$ such that $PSL(3,q)\unlhd
G_1\lhd G_2\leq G\leq Aut PSL(3,q)$ where $|G_2:G_1|$ is a prime and $G_1$ is
not line-transitive on $\spaceS$ while $G_2$ is. By  Lemma \ref{lemma:exceptional}, $(G_2, G_1, \Pi)$ is an exceptional triple and \cite[Theorem
  1.5]{gms} implies that $\Galph$ lies inside $M$ a maximal subgroup of $G$ with $M\cap L\cong PSL(3,q_0)$ where $q=q_0^a, a>3$. 

Now consider $M$ acting on $M/\Galph.$ Since $(G_2, G_1, \Pi)$ is an exceptional triple we can find an exceptional triple $(M_2,M_1,M/\Galph)$ where $PSL(3,q_0)\unlhd M_1\lhd M_2\leq M$ (see \cite[p.4]{gms}). Once more we apply \cite[Theorem 1.5]{gms} to find that $\Galph$ lies in a subfield subgroup of $M$. Iterating like this we are able to conclude that $\Galph\cong PSL(3,q_1)$ where $q=q_1^c$ and $c$ is an integer divisible by neither $2$ nor $3$. This situation corresponds precisely with the second possibility given in Theorem A.
\end{proof} 


In order to prove Theorem A it is now sufficient to proceed under the
assumptions of the following hypothesis. Our aim is to show that this
hypothesis leads to a contradiction. We will need to consider different
possibilities for a linear space $\spaceS$ having a significant prime dividing
$|PSL(3,q)|=q^3(q-1)^2(q+1)(q^2+q+1)/\mu$ where $\mu=(q-1,3)$. Recall that a significant prime is one which divides $(b,v-1)$.

\begin{hypothesis}
Suppose that $G=PSL(3,q)$ acts line-transitively but not
flag-transitively on a linear space $\spaceS$ which is not a
projective plane. Let $b,v,k,r$ be the parameters of the space. Let
$D$ be the subgroup of $PSL(3,q)$ as defined in the previous section. We
suppose, by Lemma \ref{lemma:invft}, that every involution of
$PSL(3,q)$ fixes a point. Finally we assume that $q>2$.
\end{hypothesis}

This hypothesis will hold for the remainder of the paper.

\subsection{Significant prime: $t|q^2+q+1, t\neq 3$}\label{subsection:easy1}

Suppose first that some $t|q^2+q+1, t\neq 3$ is a significant prime. Lemma
\ref{lemma:sigprime} implies that $G_\alpha\geq \hatt(q^2+q+1).3$ which is the normalizer
of a Sylow $t$-subgroup of $PSL(3,q)$. Now $\hatt (q^2+q+1).3$ is maximal in
$PSL(3,q)$ for $q\neq 4$ and so, in this case, $G_\alpha=\hatt (q^2+q+1).3$
This is a contradiction since then $\Galph$ doesn't contain any
involution, contradicting our Hypothesis.

When $q=4$ the only other possibility is that $\Galph=PSL(2,7)$ and
$v=120$. Then $17|v-1$ and by Lemma \ref{lemma:largeprime}, $k\geq 17$
which contradicts Fisher's inequality (Lemma \ref{lemma:arith}).

\subsection{Significant prime: $t=p$}\label{subsection:easy2}

Suppose now that $p$ is significant. Lemma \ref{lemma:sigprime} implies that
$G_\alpha\geq \hatt [q^3]:(q-1)^2,$ a Borel subgroup, which is the normalizer
of a Sylow $p$-subgroup of $PSL(3,q)$. Then $\Galph$ is either a Borel
subgroup or a parabolic subgroup of $PSL(3,q)$.

In the latter case the action of $G$ is $2$-transitive on points and
hence flag-transitive. Thus this case is already covered.

When $\Galph$ is a Borel subgroup $v=(q^2+q+1)(q+1)$ and, by Corollary
\ref{corollary:supergroups}, $b$ divides into
$\frac{1}{2}q(q+1)(q^2+q+1)$. This implies that $r>k>q+1$. Then
$r=\frac{v-1}{k-1}<q^2+q+1$.

Consider the set of lines through the point $\alpha$. These
lines contain all points of $\spaceS$ and so the points of
$\spaceS\backslash\{\alpha\}$ can be thought of as making up a
rectangle with dimensions $r$ by $k-1$. The area of this rectangle
(that is, the number of points in the rectangle) is
$v-1=r(k-1)=q^3+2q^2+2q$.

Now $\Galph$ has five orbits on $\spaceS\backslash\{\alpha\}$ of size $q, q, q^2, q^2$ and
$q^3$. Each of these orbits forms a rectangle of points in
$\spaceS\backslash\{\alpha\}$. Thus we have a rectangle of area
$q^3+2q^2+2q$ made out of rectangles of area $q, q, q^2, q^2$ and
$q^3$ with integer dimensions. We investigate this situation.

Write $[q^a]_l$ and $[q^a]_w$ for the length and width of the rectangles of
area $q^a$. Define the length of the rectangle of area $q^a$ to be the
dimension in the direction of the side of length $r$ in the big rectangle.
Observe first that $q \leq k-1 < [q^3]_w \leq [q^3]_l < r < q^2.$

Now if, for $a=1,2,$ there exists a rectangle such that $[q^a]_w> [q^3]_w$ or
$[q^a]_l> [q^3]_l$ then $v-1\geq pq^6$ which is a contradiction. 

Suppose that there is a rectangle of area $q^a$ such that $[q^a]_w< [q^3]_w.$
This rectangle must combine with others to make the total width $k-1$. It
either combines with a rectangle of width at least $[q^3]_w$ or it combines
with several of width less than $[q^3]_w$. Given that there are only five
rectangles in total the latter possibility can only occur for $p=2, 3$. In fact
a cursory examination can rule out the case where $p=3$. A
similar argument works if we consider lengths instead widths.

Thus, for $p>2$, we are reduced to factorising $q^2+2q+2$ in $\mathbb{Z}[x]$. But
this polynomial is clearly irreducible by Eisenstein's criterion.

If $p=2$ then a slight modification of the above argument reduces to the
factorisation problem once again and the possibility that $\Galph$ is a Borel
subgroup is excluded.

\begin{remark}
Note that we have excluded the possibility that $\Galph$ is a
parabolic or a Borel subgroup, no matter what prime is significant.
\end{remark}

\subsection{Remaining Cases}\label{subsection:table}

We wish to enumerate the remaining cases that we need to examine. First of all
note that when $q$ is small applications of Lemmas \ref{lemma:largeprime} and
\ref{lemma:arith} can be used to exclude all possibilities. We will assume from
here on therefore that $q\geq 8$. 

In addition one case in particular is worth mentioning now: When $q$ is odd and when
both $2$ and $3|(q-1)$ are significant primes.

The only maximal subgroups which have index not divisible by 2 and 3
in this case are those of type 2 and 4. Suppose that $\Galph$ lies in a
subgroup $M$ of type 2. Without loss of generality the diagonal subgroup normalized by the group
of permutation matrices isomorphic to $S_3$. Now $D$ normalizes a
Sylow $2$-subgroup of $M$. In addition $Q\in Syl_3 G$ is conjugate to $H:C_3$ where $H$ is a
diagonal subgroup, $C_3$ a group of permutation matrices. $Q$
does not normalize $D$ hence $\Galph$ contains at least two conjugates
of $D$. Since these intersect trivially, by Lemma
\ref{lemma:dconjintersect}, these generate a subgroup of index
dividing $\mu$
in the diagonal subgroup. Our group $\Galph$ must therefore be the
full subgroup of type 2.

If $\Galph$ is contained in a subgroup, $M$, of type 4 then in order
to contain an element of order $\frac{q-1}{\mu}$, $M=PSL(3,q_0), q=(q_0)^2.$
But then the index of $M$ in $G$ is even which is a contradiction.\label{psl3NINE}

Thus the cases which we need to examine are, for $q\geq 8$:
\begin{center}
\begin{tabular}{|c|c|c|}
\hline
& {\bf Significant primes } $t$ & {\bf Possible stabilizers} \\
\hline
I & $\exists t\big|(q+1), t\neq 2$ & $\hatt (q^2-1).2\leq\Galph<[q^2]:\hatt GL(2,q)$ \\
\hline
II & $\exists t\big|(q-1), t\neq 2, 3$ OR & $\Galph=\hatt(q-1)^2:S_3$ \\
 & $2, 3\big|(q-1)$ both significant & \\
\hline
III & $3\big|(q-1)$ is uniquely & $\Galph$ is a subgroup of a  \\
		 & significant  & maximal subgroup of type 2, 4, 5 or 8\\
\hline		
IV & $2\big|(q-1)$ is uniquely & $\Galph$ is a subgroup of a \\
		 & significant  & maximal subgroup of types 1, 2 or 4\\
\hline
\end{tabular}
\end{center}\label{psl3TEN}

\section{Case I: $\exists t|(q+1), t\neq 2$ significant}\label{section:qplusonesig}

In this case $\Galph$ contains a subgroup $H$ of order $2(q^2-1)/\mu$
which itself has a cyclic subgroup of size $(q^2-1)/\mu$ and $\Galph$
lies inside a parabolic subgroup of $G$.

Now observe that $H$ lies inside a copy of $\hatt GL(2,q)$ and that
$\hatt GL(2,q)$ normalizes an elementary abelian subgroup, $U$, of
$PSL(3,q)$, of order $q^2$. In its conjugation action on the
non-identity elements of $U$, $\hatt GL(2,q)$ has stabilizers of order
$q(q-1)$. Thus our group $H$ must, if it normalizes any subgroup of
$U$, normalize a subgroup of order $1+x(q+1)$ for some integer
$x$. Now for such a value to divide $q^2$, as required, $x$ must be
$0$ or $q-1$.\label{psl3ELEVEN}

Thus $\Galph=\hatt A.B$ where $A$ is trivial or of size $q^2$ and
 $H\leq B\leq GL(2,q)$. Now, in the characteristic 2 case,
 $GL(2,q)=PSL(2,q)\times (q-1)$ and $H=D_{2(q+1)}\times (q-1)$. Since
 $D_{2(q+1)}$ is maximal in $PSL(2,q)$ for all even $q\geq 8$, we know
 that $B=H$ or $B=GL(2,q)$. In the odd characteristic case,
 $GL(2,q)=/\langle-I\rangle.(PSL(2,q)\times (\frac{q-1}{2})).2$ and $H=/\langle-I\rangle.(H\times
 (\frac{q-1}{2})).2$. Now, for all odd $q>9$, $D_{2(q+1)}$ is maximal
 in $GL(2,q)$ and, once again we conclude that $B=H$ or $B=GL(2,q)$.

We need to consider the case where $q=9$ and $H<B<GL(2,q)$. In fact
this case cannot occur since the only proper subgroup of $PSL(2,9)$
containing $D_{10}$ is $A_5$, but $/\langle-I\rangle.(A_5\times (\frac{q-1}{2}))$ is
not normalized by any element of $GL(2,q)$ of non-square determinant.\label{psl3TWELVE}

Thus we can summarize the cases that we need to examine:
\begin{enumerate}
\item	$\Galph=U.\hatt(q^2-1).2$ where $U=[q^2]$;
\item	$\Galph=\hatt GL(2,q)$;
\item	$\Galph=\hatt((q^2-1).2)$.
\end{enumerate}

Note that we exclude the case where $\Galph=\hatt U:GL(2,q)$, as
then $\Galph$ is maximal parabolic and this case is already
excluded. We will consider the remaining cases in turn.

\begin{remark}
These cases also arise when $2\big|(q+1)$ is the only
significant prime (see Section \ref{section:twosig}). The arguments
given below are general and apply in that situation as well.
\end{remark}

\subsection{Case 1: $\Galph = U.(\hatt(q^2-1).2)$.}\label{subsection:qplusoneA}

Now we know that $v=\frac{1}{2}(q^2+q+1)q(q-1)$ and, since $\Galph$
lies inside a parabolic subgroup, we can appeal to Corollary
\ref{corollary:supergroups} to observe that
$$b\big| \frac{1}{8}(q^2+q+1)q(q-1)(q+1)(q-2) \ \ {\rm and} \ \ b\big| \frac{1}{4}(q^2+q+1)q(q-1)(q+1)q.$$

Thus $b\big|\frac{1}{4(2,q-1)}(q^2+q+1)q(q-1)(q+1)$ and so $4(2,q-1)q^2(q-1)/\mu$
divides $\ordGL$. For $q>7$ this means that $\GL$ lies in a parabolic
subgroup. Observe that we can presume that $U:D$ lies in $\GL$ for
some $\linel$ since $U:D$ lies in $\Galph$ and is normalized by the full parabolic
subgroup (Lemma \ref{lemma:normalizersub}).\label{psl3THIRTEEN}

Suppose that $U$ is non-normal in $\GL=\hatt A.B$ where $A$ is an elementary
abelian $p$-group and $B\leq GL(2,q)$. Then $\GL$ must lie in a parabolic subgroup which
is not conjugate to $N_G(U)$ and $|U\cap A|=q$.
If $A\backslash U$ is non-empty then $U$ acts by conjugation on these
elements with an orbit, $\Omega$, of size $q$. Then $U\cap A$ and
$\Omega$ lie inside $A$ and generate $q^2$ elements. Hence we must
have $A$ of size $q$ or $q^2$. The latter would make $\ordGL\geq
4q^3(q-1)/\mu$ which is larger than $\ordGalph$ which is a
contradiction. Hence we conclude that $|A|=q$.\label{psl3FOURTEEN}

Since $A$ is normal in $\GL$ we must have $\GL$ a subgroup of a Borel
subgroup. However in this case $U$ is normal in $\GL$. This is a contradiction.

Hence we have $U$ normal in $\GL$. Furthermore there
are no other $G$-conjugates of $U$ in $\GL$, since $U\cap U^g$ is
trivial for all $g$ in $G\backslash N_G(U)$. Hence we may appeal to
Lemma \ref{lemma:fixedset}. Then either $U:\hatt GL(2,q)$ acts
line-transitively on the fixed set of $U$, which is itself a linear
space, or this fixed set lies completely in one line. In the first
case, such an action of $U:\hatt GL(2,q)$ has a kernel $U:\hatt D$ and
corresponds to a line-transitive action of $PGL(2,q)$ with stabilizer
a dihedral group $D_{2(q+1)}$. 

Examining the results of line-transitive and 2-transitive actions of
$PGL(2,q)$ we find that there is one such action to consider. We have $q$
even and $PGL(2,q)$ acts line-transitively upon a 
Witt-Bose-Shrikhande space with line-stabilizer an
elementary abelian group of order $q$. In $PSL(3,q)$ this corresponds
to $\GL$ having order $\frac{q^3(q-1)}{\mu}$ and $b=(q-1)(q+1)(q^2+q+1)$. Then we must have,
\begin{eqnarray*}
&&k(k-1) = \frac{v(v-1)}{b} = \frac{1}{4}q(q^3-q^2+q-2) \\
&\implies& 2k(2k-2) = q^4-q^3+q^2-2q.
\end{eqnarray*}
Now observe that,
$$(q^2-\frac{1}{2}q+1)(q^2-\frac{1}{2}q-1) < q^4-q^3+q^2-2q <
(q^2-\frac{1}{2}q+2)(q^2-\frac{1}{2}q).$$ \label{psl3FIFTEEN}
Thus this case is excluded.

We can assume therefore that the set of fixed points of $U$ lies
completely in one line. This fixed set has size $\frac{1}{2}q(q-1)$
and thus $k$ is at least this large. Now the subgroups conjugate to
$U$ intersect trivially. Thus $U$ lying in $\GL$ has orbits on the
points of $\linel$ of size 1 ($\frac{1}{2}q(q-1)$ such) or $q^2$ (for
$q$ odd) or $\frac{q^2}{2}$ (for $q$ even.)  

If $k\geq q^2+\frac{1}{2}q(q-1)$ then $k(k-1)>v$ which is a
contradiction. If $k=\frac{1}{2}q(q-1)$ then $k-1=\frac{1}{2}(q+1)(q-2)$ divides into
$v-1=\frac{1}{2}(q+1)(q^3-q^2+q-2)$. This is possible only for $q\leq 4$
which is a contradiction. Thus we are left with the possibility that
$q$ is even and $k=\frac{1}{2}q(q-2)$. Once again $k-1$ dividing into
$v-1$ implies that $q\leq 4$.
\label{psl3SIXTEEN, psl3A}

\subsection{Case 2: $\Galph = \hatt GL(2,q)$}\label{subsection:qplusoneB}

Since $v=q^2(q^2+q+1)$ and $\Galph$ lies inside a parabolic subgroup, we can appeal to Corollary \ref{corollary:supergroups} to observe that
$$b\big| \frac{1}{2}q^2(q^2+q+1)(q-1)(q+1) \ \ {\rm and} \ \ b\big| \frac{1}{2}q^2(q^2+q+1)(q+1)q.$$
Thus $b\big|\frac{1}{2}q^2(q^2+q+1)(q+1)$ and so $2q(q-1)^2/\mu$ divides $\ordGL$.

This implies that, for $q>7$, $\GL$ lies in a parabolic subgroup or
$q=16$. When $q=16$ we find that the prime 4111 divides into v-1=69888
which, using Lemma \ref{lemma:largeprime}, contradicts Lemma
\ref{lemma:arith}.

Thus $\GL$ lies in a parabolic subgroup and we write $\GL=\hatt A.B$ as
usual. If $A=\{1\}$ then we must have $\GL=\hatt B\leq\hatt GL(2,q).$
Examining the subgroups of $GL(2,q)$ given in Theorem
\ref{theorem:gl2subgroups} we find that $\ordGL$ is divisible by
$\frac{|GL(2,q)|}{2\mu}$. Now if $\mu=3$ and $3$ is significant then
$\Galph$ does not lie in a parabolic subgroup. Hence we must have
$\ordGL=\frac{1}{2}\ordGalph$ with $2$ uniquely significant. But then Lemma
\ref{lemma:bdoublev} implies that any prime dividing into $v$ must
be equivalent to $1(4)$. Now in our current situation any significant
prime divides into $\frac{q+1}{2}$ thus $2$ is not a significant
prime; this is a contradiction.

Now if $1\neq g\in A$ then $|C_{PSL(3,q)}(g)|=q^3(q-1)/\mu$. Thus $B$ must
act on the non-trivial elements of $A$ with orbits of size divisible
by $q-1$. Thus $|A|=q$ or $q^2$.

If $|A|=q^2$ then $\ordGL\geq 2q^2(q-1)^2/\mu>\ordGalph$ which
cannot happen. If $|A|=q$ then $p=2$ (since, if $p$ is odd, $B$ must
act on the non-trivial elements of $A$ with orbits of size divisible
by $2(q-1)$.) For $q>4$ we must have $B$ either maximal in $GL(2,q)$
of type 4 or a subgroup of the Borel subgroup
of $GL(2,q)$. In the first case $\hatt B$ has orbits of size at least
$2(q-1)$ on the non-identity elements of $A$, thus this
case can be excluded.\label{psl3EIGHTEEN}

If $B$ lies inside a Borel subgroup of $GL(2,q)$ then $B=B_1.B_2$ where $2<B_1$ and $B_2=(q-1)^2$. In fact we must
have $|B| = \frac{q(q-1)^2}{\mu}$ since $B_2$ acts by conjugation on
the non-identity elements of $B_1$ with orbits of size $q-1$. Hence $\ordGL =
\frac{q^2(q-1)^2}{\mu}$ and $b=q(q+1)(q^2+q+1)$. Hence we must have
$$k(k-1) = q^4+q^2-q.$$\label{psl3NINETEEN}
Now observe that,
$$q^2(q^2-1)<q^4+q^2-q< (q^2+1)q^2.$$
Thus this case is excluded.

\subsection{Case 3: $\Galph=\hatt (q^2-1).2$}\label{subsection:qplusoneC}

Since $v=\frac{1}{2}q^3(q^2+q+1)(q-1)$ and $\Galph$ lies inside
a parabolic subgroup, we can appeal to Corollary
\ref{corollary:supergroups} to observe that $b$ divides into both
$$\frac{1}{4}q^3(q^2+q+1)(q-1)(q+1)q \ \ {\rm and} \ \ \frac{1}{8}q^3(q^2+q+1)(q-1)(q+1)(q^3-2q^2+2q-2).$$
Thus $b\big|\frac{1}{4(2,q-1)}q^3(q^2+q+1)(q-1)(q+1)$ and so $4(2,q-1)(q-1)/\mu$ divides $\ordGL$.

To begin with note that all cases where $11< q\leq16$ and
$q=9,19,25,31,37,64$ can be ruled out using Lemma
\ref{lemma:largeprime}. When $q=11$, Lemma \ref{lemma:largeprime}
leaves one possibility, namely that $k=444$. But then $b$ is not an
integer and so this situation can be excluded. When $q=8$, Lemma
\ref{lemma:largeprime} leaves one possibility, namely that
$k=171$. But then $k-1$ does not divide into $v-1$ and so this
situation too can be excluded.\label{psl319a} 

Using these facts, and recalling that $4(2,q-1)(q-1)/\mu$ divides
$\ordGL<\ordGalph,$ we can exclude the possibility that $\GL$ lies in
a subgroup of $PSL(3,q)$ of type 3-10. Hence we assume that
$q\geq17$ and $\GL$ lies inside a subgroup of type 1 or 2 for the rest of
this section. \label{psl320, psl330, psl334}

Now $D<\Galph$ and, by Lemma
\ref{lemma:normalizersub}, $D$ lies in $\GL$ for some line
$\linel$. We refer to Lemma \ref{lemma:fixedset} to split our
investigation into three cases:

\begin{itemize}
\item	{\bf Case 3.A:} All $G$-conjugates of $D$ in $\GL$ are
  $\GL$-conjugate and the fixed set of $D$ is a linear-space acted on
line-transitively by $\hatt GL(2,q)$, the normalizer of $D$.

\item	{\bf Case 3.B:} All $G$-conjugates of $D$ in $\GL$
are $\GL$-conjugate and the fixed points of $D$, of which there are
$\frac{1}{2}q(q-1)$, lie on one line; 

\item	{\bf Case 3.C:} $\GL$ contains at least two $\GL$-conjugacy
  classes of $G$-conjugates of $D$.
\end{itemize}

\subsubsection{Case 3.A}

This situation corresponds to a line-transitive action of
$PGL(2,q)$ with stabilizer $D_{2(q+1)}$. Then Theorem
\ref{theorem:psl2} implies that
$p=2$ and the fixed set of $D$ is a Witt-Bose-Shrikhande
space. The corresponding line-stabilizer in $PGL(2,q)$ has size $q$
and so $\ordGL$ is divisible by $\frac{q(q-1)}{\mu}$ in
$PSL(3,q)$. Suppose that $\ordGL=\frac{q(q-1)}{\mu}$ and so
\begin{eqnarray*}
&&k(k-1)=\frac{v(v-1)}{b}=\frac{1}{4}(q^6-q^5+q^4-2q^3+2q^2-2q) \\
&\implies&(2k)(2k-2)=q^6-q^5+q^4-2q^3+2q^2-2q.
\end{eqnarray*}
But now observe that
$$(q^3-\frac{1}{2}q^2+\frac{3}{8}q+2)(q^3-\frac{1}{2}q^2+\frac{3}{8}q)< 2k(2k-2) < (q^3-\frac{1}{2}q^2+\frac{3}{8}q)(q^3-\frac{1}{2}q^2+\frac{3}{8}q-2).$$
For $q>16$ this gives a contradiction.\label{psl3TWENTY-ONE}

The only other possibility is that $\ordGL=\frac{2q(q-1)}{\mu}$ and
$[q]\times\frac{q-1}{\mu}=\GL\cap C_G(D)$. This implies that $\GL$
lies inside a parabolic subgroup of $PSL(3,q)$. 

Now $[q]\times\frac{q-1}{\mu}$ is normal in $\GL$ and so $[q]$ is
normal in $\GL$ and $\GL$ lies inside a Borel subgroup of
$PSL(3,q)$. Then $D$ acts on the normal subgroup of $\GL$ of order
$2q$. Furthermore $D$ centralizes at most $q$ of these elements and
has orbits on the rest of
size at least $\frac{q-1}{\mu}$. These orbits intersect cosets of
$[q]\unlhd C_G(D)\cap\GL$ with a size of at most $1$. This gives a
contradiction.\label{psl3TWENTY-TWO}

\subsubsection{Case 3.B}

Observe that all
$PSL(3,q)$-conjugates of $D$ intersect trivially. Observe too that all
elements of $\Galph$ are of form $TS$ where $T\in \hatt (q^2-1)$ and
$S^2$ lies in $D$. Then $(TS)^2$ lies in $D$ and hence if $E$ is some
other conjugate of $D$ then $E\cap \Galph$ is of size at most
$(2,q-1)$. Thus the orbits of $D$ on $\linel$, a line which it fixes, are
either of size $\frac{q-1}{(2,q-1)\mu}$ or of size $1$ and there are
$\frac{1}{2}q(q-1)$ of these. We conclude that $k$ is a multiple of
$\frac{q-1}{(2,q-1)\mu}$.

Now we find that $(v-1,|G|)=\frac{q+1}{(2,q-1)}$. Since
$\frac{q-1}{(2,q-1)\mu}\big|k$ and
$b=\frac{v(v-1)}{k(k-1)}$ divides into $|G|$ then
$b\big|\frac{\mu}{2}(q^2+q+1)q^3(q+1)$. 

Thus, for $q\not\equiv 1(3)$, $\ordGL=2(q-1)^2\geq 512$. If
$q\equiv 1(3)$ then $\ordGL=\frac{2}{9}(q-1)^2.a$ where $a=1,2$ or
$3$.

Suppose first that $p$ is odd. Consider the possibility that
$\GL$ lies inside a subgroup of type 2 and not 
in a parabolic subgroup. So $\GL$ is a subgroup of $\hatt (q-1)^2:S_3$
and must have either $3$ or $S_3$ on top. The former case is
impossible as then $b$ does not divide into
$\frac{\mu}{2}(q^2+q+1)q^3(q+1)$. Now $\GL=\hatt (A\times A):S_3$ or
$(\frac{A}{\mu}\times\frac{A}{\mu}):S_3$. Then, since $\GL$ must
contain a subgroup conjugate to $D$, we find that
$\GL=(\frac{q-1}{\mu}\times \frac{q-1}{\mu}):S_3$, $\mu=3$ or
$\GL=\hatt(q-1)^2:S_3$. The latter case violates Fisher's
inequality and can be excluded. In the former case $\GL$ contains at
most $q+2$ involutions. Appealing to Lemma \ref{lemma:invineq}, we
observe that
$$k\leq \frac{r_gv}{n_g}+1=\frac{1}{2}q(q+2)(q-1)+1.$$
This means that $b=\frac{v(v-1)}{k(k-1)}>q^5(q-3)$ which is a contradiction.\label{psl3THIRTY-ONE}

Thus $\GL$ lies inside a parabolic subgroup; in fact $\GL$ is
isomorphic to a subgroup of $\hatt GL(2,q)$. In order for Fisher's
inequality to hold, we must have one of the following cases:
\begin{itemize}
\item	$b=\frac{1}{2}q^3(q^2+q+1)(q+1)$ and
  $|\GL|=\frac{2(q-1)^2}{\mu}$. Thus $\GL$ is isomorphic to a
subgroup of $\hatt GL(2,q)$ of type 4 (in which case $\GL$ contains
more than one $\GL$-conjugacy class of $G$-conjugates of $D$ which is a contradiction) or $\GL$ is
isomorphic to a subgroup of type 6 or 7. This latter case requires
that $2(q-1)$ divides into $24$ or $60$. These possibilities have
already been excluded.\label{psl3TWENTY-FIVE}

\item	$b=\frac{3}{4}q^3(q^2+q+1)(q+1)$. Hence $\ordGL=\frac{4}{9}(q-1)^2$ and $q\equiv 7(12)$. Thus $\GL$ is
isomorphic to a subgroup of type 6 or 7 in $\hatt GL(2,q)$ and
$\frac{4(q-1)}{3}$ must divide $24$ or $60$. This is
impossible.\label{psl3TWENTY-SIX} 

\item	$b=\frac{3}{2}q^3(q^2+q+1)(q+1)$. Then
  $\ordGL=\frac{2}{9}(q-1)^2$ and $q\equiv 1(3)$. Thus $\GL$ is isomorphic to a
subgroup of $\hatt GL(2,q)$ of type 4, 6 or 7. \label{psl3TWENTY-SEVEN}

If $\GL$ is isomorphic to a subgroup of $\hatt GL(2,q)$ of type 4
then $r_g\leq\frac{q+8}{3}$. Using Lemma \ref{lemma:invineq} we see
that\label{psl3TWENTY-EIGHT} 
$$k\geq \frac{n_g(v-1)}{br_g}> q^2(q-9).$$
Since $(k-1)^2<v$ this implies that
$$q^4(q-9)^2<\frac{1}{2}q^3(q^2+q+1)(q-1)$$
which means that $q<31$. Then $q=25,$ but this possibility
has already been excluded using Lemma
\ref{lemma:largeprime}.\label{psl3TWENTY-NINE} 

If $\GL$ is isomorphic to a subgroup of $\hatt GL(2,q)$ of type 6 or 7
then we require that $\frac{2(q-1)}{3}$ divides into $24$ or
$60$. Hence $q=31$ or $37$. These possibilities have already been
excluded.\label{psl3THIRTYB}
\end{itemize}

If $p=2$ then, in order for Fisher's inequality to hold and so that
$4(q-1)/\mu$ divides into $\ordGL$, we have $\ordGL =
\frac{4}{9}(q-1)^2$ and $q\equiv 1(3)$. Thus $\GL$ lies inside a
parabolic subgroup of $PSL(3,q)$ and $\GL=\hatt A.B$ as usual.

If $A$ is trivial then $\GL$ is a subgroup of type 2 in $\hatt
GL(2,q)$. Then $\GL$ has a normal $2$-group and, by Schur-Zassenhaus,
$\GL$ also contains a subgroup of size
$\frac{(q-1)^2}{9}$. This subgroup has orbits in its conjugation
action on $2$-elements of $\GL$ of size at least $\frac{q-1}{3}$. This
implies that $\ordGL$ is divisible by $\frac{q(q-1)^2}{9}$ which is a
contradiction. \label{psl3THIRTY-THREE}

If $A$ is non-trivial then $\GL$ must have orbits in its conjugation
action on non-identity elements of $A$ of size at least
$\frac{q-1}{3}$. Once again this implies that $\ordGL$ is divisible by
$\frac{q(q-1)^2}{9}$ which is a contradiction.

\subsubsection{Case 3.C}
Now consider the possibility that $\GL$ contains at least two
$\GL$-conjugacy classes of $G$-conjugates of $D$.

Suppose first that $\GL$ is a subgroup of
$\hatt(q-1)^2:S_3$ and does not lie in a parabolic subgroup. We know
that $q$ is odd since $4(2,q-1)(q-1)/\mu$ divides into $\ordGL$. Since
$\GL$ is not in a parabolic subgroup we must have a non-trivial part of $S_3$ on
top, of order $3$ or $6$. Thus all $G$-conjugates of $D$ in $\GL$ are
$\GL$-conjugate which is a contradiction. 


Thus we may conclude that $\GL$ is in a parabolic subgroup. Write
$\GL=\hatt A.B$ as usual.  If $A$ is trivial then, referring to Lemma
\ref{lemma:dconjugates}, we conclude that $\GL$ is a subgroup of
$\hatt GL(2,q)$ of type 2,4 or 5. If $\GL$ is of type 5 then $q=49$
and this can be ruled out using Lemma \ref{lemma:largeprime}. \label{psl3THIRTY-FIVE} 

If $\GL$ is of type 2 and not of type 4 then it must contain
non-trivial $p$-elements. Some conjugate of $D$ in $\GL$ must have orbits in
its conjugation action on these elements of size
$\frac{q-1}{\mu}$. Thus $A_1:\frac{q-1}{\mu}\leq\ordGL$ where $A_1$ is
a $p$-group of size divisible by $q$. We will
consider this possibility together with the case when $A$ is non-trivial.

So suppose that $A$ is non-trivial. Now either all $G$-conjugates of $D$
in $\GL$ lie in $C_G(A)$ or else $|A|\geq
q$. Consider the first possibility. In this case $A:D$ and $A:E$ lie
inside $C_G(A)$ where $E$ is a $G$-conjugate of $D$. Now $C_G(A)\leq
C_G(g)$ for $g$ an element or order $p$. Since $C_G(g)\cong
[q^3]:\frac{q-1}{\mu}$, we know that $D$ and $E$ are conjugate in
$C_G(A)\cap\GL$ by Schur-Zassenhaus. This is a contradiction and so we
assume that $|A|\geq q;$ thus, in both cases that we have considered
so far, $Q:D\leq \GL$ where $Q$ is a $p$-group of order divisible
by $q$.

Now let $E$ be a $G$-conjugate of $D$ in $\GL$ which is not
$\GL$-conjugate to $D$. Suppose $E\cap (Q:D)$ is non-trivial
and $1\neq h\in E\cap (Q:D)$. Then $h$ lies inside a
$Q:D$-conjugate of $D$ by applying Sylow theorems to $Q:D$. But
this is impossible since Lemma \ref{lemma:dconjintersect} implies that
either $E=D$ or $E\cap D$ is trivial. Hence $\ordGL\geq
\frac{q(q-1)^2}{\mu^2}>\ordGalph$ which is also
impossible.\label{psl3THIRTY-SIX}

Finally we must consider the possibility that $\GL$ is of type 2 in
$\hatt GL(2,q)$; that is, $\GL$ is a subgroup of
$\hatt(q-1)^2:2$. We must have $q$ odd since $4(2,q-1)(q-1)/\mu$
divides into $\ordGL$. Furthermore the $G$-conjugates of $D$ in
$\hatt(q-1)^2:2$ normalize each other and so $\frac{(q-1)^2}{\mu^2}$
divides into $\ordGL$. There are three possibilities to consider:
\begin{itemize}
\item	$\GL\leq \hatt (q-1)^2$. In this case $\GL$ contains at most $3$ involutions. Appealing to Lemma \ref{lemma:invineq}, we observe that
$$k\leq \frac{r_gv}{n_g}+1=\frac{3}{2}q(q-1)+1.$$
This is too small to satisfy $b=\frac{v(v-1)}{k(k-1)}$ hence we have a contradiction.
\item $\GL = (\frac{q-1}{\mu}\times \frac{q-1}{\mu}):2$. Then $\GL$ contains $\frac{q+8}{3}$ involutions. Once again using Lemma \ref{lemma:invineq}, we observe that
$$k\leq \frac{r_gv}{n_g}+1=\frac{1}{6}q(q+8)(q-1)+1.$$
But this is too small to satisfy $b=\frac{v(v-1)}{k(k-1)}$ hence we have a contradiction.\label{psl3THIRTY-SEVEN}
\item $\GL = \hatt ((q-1)\times (q-1)):2$. Then $\GL$ contains $q+2$ involutions and we have that,
$$k\leq \frac{r_gv}{n_g}+1=\frac{1}{2}q(q+2)(q-1)+1.$$
Once again this is too small to satisfy $b=\frac{v(v-1)}{k(k-1)}$.\label{psl3THIRTY-ONEb}
\end{itemize}

Hence we may conclude that no line-transitive actions exist with primes dividing $q+1$ significant.

\section{$\Galph=\hatt (q-1)^2:S_3$}

In this case $v=\frac{1}{6}q^3(q+1)(q^2+q+1)$ and any significant prime
$t$ must divide into $q-1$.

\begin{remark}
The argument in this section deals with Case II in our analysis of
significant primes.
\end{remark}

Note first that, by using Lemma \ref{lemma:largeprime}, we can
assume that $q>25$ and that $q\neq 31$, $37$, $43$, $49$, $64$, $109$ or
$271$.\label{psl3THIRTY-EIGHT} Furthermore a conjugate of $D$ lies in
$\Galph$ and $D$ is normalized by $\hatt GL(2,q)$. Thus, by Lemma
\ref{lemma:normalizersub}, a conjugate of $D$ lies inside $\GL$. We
split into three cases:
\begin{itemize}
\item	{\bf Case A}: A $G$-conjugate of $D$ is normal in $\GL$ and $\GL$ contains
no other $G$-conjugates of $D$;
\item	{\bf Case B}: A $G$-conjugate of $D$ is normal in $\GL$ and $\GL$ contains
other $G$-conjugates of $D$. Thus $\ordGL$ is divisible by
$(\frac{q-1}{\mu})^2$ and so $b$ divides into $6\mu v$;
\item 	{\bf Case C}: All $G$-conjugates of $D$ in $\GL$ are non-normal in $\GL$.
\end{itemize}

We examine these possibilities in turn.

\subsection{Case A}
In this case we know, by Lemma \ref{lemma:fixedset}, that either
$\hatt GL(2,q)$ acts line-transitively on the linear-space which is the
fixed set of $D$ or all fixed points of $D$ lie on a single line. The
first possibility cannot occur however as this would correspond to
$PGL(2,q)$ acting line-transitively on a linear-space (possibly having
$k=2$ and so being a $2$-homogeneous action) with line-stabilizer a
dihedral group of size $2(q-1)$ which is impossible. Hence we may
assume that all fixed points of $D$ lie on a single line. There are
$\frac{1}{2}q(q+1)$ of these. 

If $E$ is some other conjugate of $D$ then $E\cap\Galph$ is of size at
most 2. We conclude that\label{psl3FORTY}
$k=\frac{1}{2}q(q+1)+n\frac{q-1}{2\mu}$ for some integer $n$. This
implies that $k-1$ is divisible by $\frac{q-1}{2\mu}$. Now, since
$v-1=\frac{q-1}{2}\frac{q^5+3q^4+5q^3+6q^2+6q+6}{3}$, we observe that\label{psl3FORTY-ONE}
$b\big| (q^5+3q^4+5q^3+6q^2+6q+6)v$. Now, for $p$ odd, $(|G|,
q^5+3q^4+5q^3+6q^2+6q+6)$ is a power of $3$, hence $3$ is the only
significant prime and $3|q-1$. For $p=2$, $(|G|,
q^5+3q^4+5q^3+6q^2+6q+6)$ is divisible, at most, by the primes $2$ and
$3$. However we know that $2$ is not a significant prime here thus,
again, $3$ is the only significant prime. Note that
$q^5+3q^4+5q^3+6q^2+6q+6$ is divisible by $27$ if and only if $q\equiv
28(81)$. Thus, if $3^a$ is the highest power of $3$ in $q-1$ then
$a\neq 3$ implies that $b|27v$. If $a=3$ then we know already that
$b|81v$.\label{psl3FORTY-TWO}

This case will be completed below.

\subsection{Case A and B}

Now we examine the remaining possibilities of Case A along with Case B. Thus $\GL<\hatt GL(2,q)$ and one of the following holds:
\begin{itemize}
\item	$q\equiv 28(81)$, $\frac{2(q-1)^2}{81}$ divides into $\ordGL$
  and $\GL$ contains precisely one $G$-conjugate of $D$;
\item	$q\equiv 1(3)$, $\frac{2(q-1)^2}{27}$ divides into $\ordGL$
  and $\GL$ contains precisely one $G$-conjugate of $D$;
\item	$\frac{(q-1)^2}{\mu^2}$ divides into $\ordGL$ and $\GL$
  contains more than one $G$-conjugate of $D$.
\end{itemize}


Observe also that $k(k-1)=\frac{v(v-1)}{b}$ is even and that 
$$|v(v-1)|_2=\frac{(q,2)}{4}|q^3(q+1)(q-1)|_2.$$
Thus if $p$ is odd then we need $\ordGL$ divisible by $8(q-1)/\mu$.

Suppose that $\GL$ is a subgroup of $\hatt GL(2,q)$ of type 6 or
7. Since $q>25,$ Lemma \ref{lemma:dconjugates} implies that $\GL$
contains at most one conjugate of $D$. Thus $\frac{2(q-1)}{9}$
must divide 24 or 60 or $\frac{2(q-1)}{27}$ divides 24 or 60 and
$q\equiv 28 (81)$. The prime powers we need to check are,
therefore, 13, 19, 31, 37, 109 and 271. These cases are already all
excluded.\label{psl3FORTY-THREE}

If $\GL$ lies inside a group of type 3 then $\GL$
contains at most one conjugate of $D$ and either $q\cong 28 (81)$
and $\frac{2(q-1)}{27}$ divides into 4 or $\frac{2(q-1)}{9}$ divides
into 4. Both yield values for $q$ which are less than 25 and so can be
excluded.\label{psl3FORTY-FOUR}

Suppose that $\GL$ is a subgroup of $\hatt GL(2,q)$ of type 5,
$\GL\cong\hatt<SL(2,q_0), V>$. Then
$\frac{(q-1)^2}{81}$ divides into $2q_0(q_0^2-1)\frac{q_0-1}{3}$ and
so $q-1$ divides into $54(q_0^2-1)$. For $q\geq q_0^3$ we find that
this is impossible for $q_0>2$. If $q_0=2$ then $q<32$ and so all
cases have been excluded. For $q=q_0$, $\ordGL<\ordGalph$ implies a
contradiction. For $q=q_0^2$, $\ordGL<\ordGalph$ implies that
$\sqrt{q}\leq 5$ and all possibilities have been
excluded.\label{psl3FORTY-FIVE}

Suppose that $\GL$ lies inside a parabolic subgroup of $\hatt GL(2,q)$ and not
of type 4. Then $\ordGL$ is divisible by $p$ for $q=p^a$, integer $a$. If
$\ordGL$ is divisible by $\frac{(q-1)^2}{\mu^2}$ then $\GL$ has orbits on
the non-identity elements of its normal $p$-Sylow subgroup divisible
by $\frac{q-1}{\mu}$. Thus $\GL$ contains the entire Sylow
$p$-subgroup of $\hatt GL(2,q)$ and $\ordGL\geq
q\frac{(q-1)^2}{\mu^2};$ this
implies that $q<6\mu$ which is impossible. So assume that $3\big|(q-1)$ is
the only significant prime. If $\frac{2(q-1)^2}{81}$
divides into $\ordGL$ we must have $p=2$ and $\GL=\hatt A:B$ where $A$ is a
non-trivial $2$-group. Then $q\geq 2^a$ and $q-1$ has a primitive prime
divisor $s$ greater than 3 and $\frac{s(q-1)}{3}$ divides into
$|B|$. Then $B$ acts on the non-identity elements of $A$ by
conjugation with orbits of size divisible by $s$ and so $|A|=q$. Thus
$\ordGL$ is divisible by $\frac{q(q-1)s}{3}$ which means $s$ must be $5$ and
so $q=16$. This is already excluded.\label{psl3FORTY-SIX}

We are left with the possibility that $\GL$ is a subgroup of $\hatt
GL(2,q)$ of type 4. If $2$ is significant then $p$ is odd and $\GL$ contains at most
3 involutions since $\GL\leq \hatt (q-1)^2$. By Lemma \ref{lemma:invineq} we know that $k\leq
\frac{3v}{n}+1=\frac{1}{2}q(q+1)+1$. This is inconsistent with our
value for $b$. \label{psl3FORTY-SEVEN}If $2$ is not significant then
$\ordGL=2|D|e$ where $e$ is a constant dividing $q-1$. Then the number
of involutions in $\GL$ is at most $e+3$. We appeal to Lemma
\ref{lemma:invineq} to conclude that,
$$k\leq \frac{r_gv}{n_g}+1=\frac{(e+3)(q+1)q}{6}+1.$$
Thus,
$$\frac{3(q-1)}{e}=\frac{b}{v}=\frac{v-1}{k(k-1)}\geq\frac{6(q^6+2q^5+2q^4+q^3-6)}{(e+3)^2q^2(q^2+3q+2)}>\frac{6q^2}{(e+3)^2}.$$
This implies that $\frac{(e+3)^2}{e}>2q$ and so $e+15>2q$. Since $e<q$
this must mean that $q<15$ which is a contradiction.\label{psl3FORTY-EIGHT}

\subsection{Case C}
Finally we consider the possibility that no conjugate of $D$ is normal
in $\GL$. We must have at least two conjugates of $D$ in $\GL$ and so $\ordGL>\frac{(q-1)^2}{\mu^2}$. 

Suppose first that $\GL$ lies in a parabolic subgroup. Then
$\GL=\hatt A.B$ where $A$ is an elementary abelian $p$-group, $B\leq
GL(2,q)$.

Suppose that $A$ is trivial and refer to Lemma
\ref{lemma:dconjugates}. Then $\GL$ lies in a subgroup of $\hatt GL(2,q)$
of types 2, 4 or 5. If $\GL$ lies in a subgroup of type 5 then
$\GL\geq SL(2,q)$ in which case $\ordGL>\ordGalph$ which is a
contradiction.

If $\GL$ lies in a subgroup of $\hatt GL(2,q)$ of type 4 then
conjugates of $D$ in $\GL$ normalize each other and so
$\frac{(q-1)^2}{\mu^2}$ divides into $\ordGL$. In this case some
conjugate of $D$ must be normal in $\GL$ which is a contradiction.

If $\GL$ lies in a subgroup of $\hatt GL(2,q)$ of type 2 then we must
have $p$ dividing $\ordGL$ otherwise all conjugates of $D$ are normal
in $\GL$. But then some conjugate of $D$ acts by conjugation on the
non-trivial elements of the normal $p$-subgroup with orbits of size
$\frac{q-1}{\mu}$. Thus $q$ divides $\ordGL$ and $\GL$ has a normal
subgroup $Q$ of size $q$. We will deal with this situation at the end
of the section.

Thus $A$ is non-trivial. Suppose that all conjugates of $D$ in $\GL$
centralize all elements of $A$. Then these conjugates lie in a
subgroup of order $q^3(q-1)/\mu$. Now if $\GL\cap C_G(A)$ only
contains $p$-elements centralized by $D$ then $\GL\cap C_G(A)$
contains only one conjugate of $D$. By our supposition this means that
$\GL$ contains only one conjugate of $D$ which is a
contradiction. Thus $\GL\cap C_G(A)$ contains $p$-elements not
centralized by $D$. Then the normal $p$-subgroup of $\GL\cap C_G(A)$
has size $|A|+n\frac{q-1}{\mu}$ for some $n$. Thus $\GL\geq Q:D$ for a
$p$-group $Q$ of size at least $q$.

If a conjugate of $D$ in $\GL$ does not act trivially in its action on
elements of $A$ then $A$ must be of order divisible by $q$. Once again
$\GL\geq Q:D$ where $|Q|\geq q$. We deal with this situation
at the end of the section.

Now suppose that $\GL$ lies inside a subgroup of $PSL(3,q)$ of type
2. In order for there to be two conjugates, $D$ and $E,$ of $D$ in
$\GL$ we must have $D,E$ in $\hatt (q-1)^2$. Hence
$\frac{(q-1)^2}{\mu^2}\big|\ordGL$. For $D, E$ to be non-normal, we
must have $\GL\geq  (\frac{q-1}{\mu}\times\frac{q-1}{\mu}):3$. If $2$
is significant then $p$
is odd and $\GL\leq \hatt(q-1)^2:3$ and $\GL$ contains at most 3
involutions. By Lemma \ref{lemma:invineq}, we know that $k\leq
\frac{3v}{n}+1=\frac{1}{2}q(q+1)+1$. This is inconsistent with our
value for $b$. If $2$ is not significant then
$\GL=(\frac{q-1}{3}\times\frac{q-1}{3}):S_3$ and $b=3v$. 

When $p$ is odd, $\GL$ contains at
most $q+2$ involutions and, by Lemma \ref{lemma:invineq}, this implies that 
$k\leq \frac{(q+2)v}{q^2(q^2+q+1)}+1$.  We therefore conclude that 
$$k(k-1)\leq\frac{q(q+1)(q+2)(q+3)(q^2+2)}{36}.$$
However this implies that $\frac{b}{v}=\frac{v-1}{k(k-1)}>4$ which is
a contradiction.

When $p=2,$ $\GL$ contains at most $q-1$ involutions and  we find that
$k(k-1)\leq\frac{1}{36}q^3(q^3+6)$. Once again $\frac{b}{v}=\frac{v-1}{k(k-1)}>4$ which is a contradiction.\label{psl3FORTY-NINE}

If $\GL$ lies inside a subgroup of $PSL(3,q)$ of type 4 or 5 then we have two
possibilities. If $\GL = A_6.2$ or $A_7$ then, in order to satisfy
$\ordGL>\frac{(q-1)^2}{\mu^2}$, we must have $q=25$. This has already
been excluded. If $\GL$ contains a subgroup of index less than or
equal to 3 isomorphic to $PSU(3,q_0)$ or $PSL(3,q_0)$ where $q=q_0^a$ then
we require that $q_0^3(q_0^2-1)(q_0^3-1)<6(q-1)^2$. Thus we need $q\geq
q_0^4$. This implies that either $\frac{q-1}{\mu}$ does not divide into
$\ordGL$ or that $q=64$. Both cases give contradictions.\label{psl3FIFTY}

If $\GL$ lies inside a subgroup of $PSL(3,q)$ of type 6,7,8 or 10 then
$\frac{(q-1)^2}{\mu^2}<360$. This implies that $q\leq 19$ or $q\equiv
1(3)$ and $q\leq 49$. All of these cases have been excluded
already.\label{psl3FIFTY-ONE}

If $\GL$ is in a group of type 9 then $\ordGL<\ordGalph$ implies that
$\GL$ is a proper subgroup. Since $\ordGL>\frac{(q-1)^2}{\mu^2}$ we
must have $\GL\leq [q]:(q-1)$. Thus $\GL=A:B$ where $A\leq [q]$, $B\leq
(q-1)$. All conjugates of $B$ in $\GL$ are $\GL$-conjugate and $B$
contains a conjugate of $D$. Thus $\frac{q-1}{\mu}$ divides into
$|B|.$ Since $B$ acts semi-regularly on the non-trivial elements of
$A$ this means that $|A|=q$. Once more we conclude that $\GL$ has a
normal subgroup of order $q$.

We have reduced all cases to the situation where $\GL\geq Q:D$ where
$Q$ is a $p$-group of order divisible by $q$. Observe that all
conjugates of $D$ in $Q:D$ are $\GL$ conjugate. If $\GL$ contains $E$,
another $G$-conjugate of $D$ which
is not $\GL$-conjugate, then ${E\cap (Q:D)}$ is trivial; hence $\ordGL\geq
\frac{q(q-1)^2}{\mu^2}$ which is too large. Thus all $G$-conjugates of
$D$ in $\GL$ are $\GL$-conjugate and we can apply Lemma
\ref{lemma:fixedset} as in Case A. As in Case A this implies that $3$
is uniquely significant and either $2\frac{(q-1)^2}{81}\big|\ordGL,
q\equiv 28(81)$ or $2\frac{(q-1)^2}{27}\big|\ordGL, q\equiv 1(3)$. If
$p$ is odd then this means that either $q<81$ and $q\equiv 28(81)$ or
$q<27$ and $q\equiv 1(3)$. If $p=2$ then this means that either
$q<162$ and $q\equiv 28(81)$ or $q<54$ and $q\equiv 1(3)$. All such
possibilities have already been excluded.

Hence we may conclude that no new line-transitive action of $PSL(3,q)$
exists where $\Galph=\hatt (q-1)^2:S_3$. 

\section{Case III: $3|q-1$ is uniquely significant}\label{section:threesig}

In this case $\Galph$ lies inside a subgroup of $PSL(3,q)$ of type 2,
4, 5 or 8. Note that $\Galph$ must contain a subgroup of type $3^2.Q_8$ so, in particular, if $\Galph$ lies inside a subgroup of type 8 then $\Galph\cong 3^2.Q_8$.

\subsection{Case 1: $\Galph$ is a proper subgroup of a group of
type 2}

Then $\Galph=A.B$ where $B=C_3$ or $S_3$
and $A=\hatt(u\times u)$ (this structure for $A$ follows since it is
normalized by $C_3$.) We can conclude, using Corollary
\ref{corollary:supergroups}, that $B=S_3$. Now observe that $A.2$ lies
inside a copy of $\hatt GL(2,q),$ hence is centralized by $Z(\hatt
GL(2,q))$. Thus, by Lemma \ref{lemma:normalizersub}, $A.2$ lies in
$\GL$. Thus $\ordGL=2|A|$ or $\ordGL=4|A|$ while $b\big|3v$. When $p=2$ we
know that $v-1$ is odd. Since $k(k-1)$ is even and
$\frac{b}{v}=\frac{v-1}{k(k-1)}$, this means that $\ordGL=4|A|$ and $b=\frac{3}{2}v$.

Consider first the case where $b=\frac{3}{2}v$. Then
$\frac{b}{v}=\frac{3}{2}=\frac{v-1}{k(k-1)}$ and so
$$k(k-1)=\frac{2}{3}(v-1)=\frac{1}{9u^2}[q^8-q^6-q^5+q^3-6u^2].$$
Now observe that, for $q>8$,
$$
[\frac{1}{3u}(q-1)(q^3+q^2+\frac{1}{2}q)+\frac{1}{2}][\frac{1}{3u}(q-1)(q^3+q^2+\frac{1}{2}q)-\frac{1}{2}]
> \frac{2}{3}(v-1); $$

$$[\frac{1}{3u}(q-1)(q^3+q^2+\frac{1}{2}q)\frac{1}{3}][\frac{1}{3u}(q-1)(q^3+q^2+\frac{1}{2}q)-\frac{2}{3}]
< \frac{2}{3}(v-1).$$

Since $\frac{1}{3u}(q-1)(q^3+q^2+\frac{1}{2}q)=\frac{1}{6}a$ for some
integer $a$, this is a contradiction. Thus $p$ is odd and $b=3v$.\label{psl3FIFTY-TWO}

Now suppose that $4$ does not divide into $u$. Then $|\Galph|_2\leq 8$
while $|G|_2\geq 16,$ hence $v-1$ is odd. This implies that
$|b|_2<|v|_2$ which is a contradiction. Hence $12|u$.

Now $\GL=\hatt(u\times u).2<\hatt (q-1)^2:2<\hatt GL(2,q)$ and so
contains at most $u+3$ involutions. We appeal to Lemma \ref{lemma:invineq} to observe that,
$$k\leq\frac{(u+3)q(q+1)(q-1)^2}{6u^2}+1.$$
We can conclude therefore that, for $u\geq 12$,
$$k(k-1)\leq \frac{q^2(q+1)^2(q-1)^4(u+3)(u+4)}{36u^4}.$$
This is strictly smaller than $\frac{v-1}{3}$ which is a
contradiction. \label{psl3FIFTY-THTREE}

\subsection{Case 2:  $\Galph$ lies inside a subgroup of type 4 or 5}

We refer to Lemma \ref{lemma:psl3subgroups}. We cover the case where $\Galph=3^2.Q_8$ below. Clearly $3^2.Q_8$ is not a subgroup of $A_6$ so $\Galph$ is not isomorphic to $A_6$. Consider the possibility that $\Galph$ is isomorphic to $A_6.2$ or $A_7$ and $p=5$. We exclude $q=25$ using Lemma \ref{lemma:largeprime}.\label{psl3FIFTY-THTREEa}

Observe that, since $3$ divides $q-1$, there is a group of order $3$ normal in a group isomorphic to $\hatt(q-1)^2$. Hence a line-stabilizer contains a subgroup of order $3$ or else contains the  group $\hatt(q-1)^2$ (by Lemma \ref{lemma:normalizersub}). The latter possibility is not possible, hence we may assume that $3\big|\ordGL$. We may therefore conclude that $b=3v$ or $b=\frac{3}{2}v$. 

Now suppose that $m$ is an integer dividing $v$ and $b=\frac{3}{x}v$ where $x$ is 1 or 2. We have that
\begin{eqnarray*}
&&\frac{v-1}{k(k-1)}=\frac{3}{x} \\
&\implies& 3k(k-1)+x\equiv 0 \ (mod \ m) \\
&\implies& 36k^2 -36k +12x\equiv 0 \ (mod \ m) \\
&\implies& 9(2k-1)^2\equiv 9-12x  \ (mod \ m) \\
\end{eqnarray*}
Thus $9-12x$ is a square modulo $m$ and $m$ is not divisible by 3. If
$\Galph=A_6.2$ then we know that $25$ divides $v$. For both values of
$x$ we find that $9-12x$ is not a square modulo $25$.

Thus we assume that either
$\Galph=PSL(3,q_0), q=q_0^a, 3\big|q_0-1, a\not\equiv 0 (mod \ 3);$ or 
$\Galph=PSU(3,q_0), q=q_0^a, 3\big|q_0+1, a\not\equiv 0 (mod \ 6).$

Then in the first instance we have a subgroup of $\Galph$,
$\hatt(q_0-1)^2$; in the second instance we have a subgroup of $\Galph$,
$\hatt(q_0+1)^2$. Such subgroups are normal in the subgroup of
$PSL(3,q)$, $\hatt (q-1)^2$. Thus these subgroups of $\Galph$ lie in
$\GL$ and we may conclude that $b\big|3v$. Once again when $p=2$ we
know that $v-1$ is odd and so $b=\frac{3}{2}v$.


We know that $q_0^3\big|\ordGL$, hence $\GL$ is not a
subgroup of a group of type 2,3,6,7,8 or 10. If $\GL$ is a subgroup of
a group of type 9 then $\frac{(q_0^3\pm1)}{3}\big|(q^2-1)$. Since
$q=q_0^a, a\not\equiv 0 (3)$ we must have $q_0=2$ and
$\Galph=PSU(3,2)$\label{psl3FIFTY-SEVENa}. But then $\ordGalph = 72$
which is the same size as in Case 1 with $u=6$. The arguments
given there exclude both $b=3v$ and $b=\frac{3}{2}v$.

If $\GL$ is only a subgroup of a group of type 4 or 5 then either
$\GL=A_6.2$ or $A_7$ (and $25$ divides into $v$ which is a
contradiction), or $\GL$ is one of $PSL(3,q_1)$ or $PSU(3,q_1)$. Since
$b|3v$ we must have $q_0=q_1$ and $\frac{q_0^3+1}{q_0^3-1}$ equal to $3$ or
$\frac{3}{2}$. This is impossible.

Thus $\GL$ is a subgroup of a parabolic subgroup. Then we require that
$(q_0^3\pm 1)\big|(q^2-1)(q-1)$. This implies that $q_0=2$ which can be
excluded as in Case 1 setting $u$ to be $6$.\label{psl3FIFTY-EIGHTa}

\subsection{Case 3: $\Galph\cong 3^2.Q_8$}

Note that $p$ is odd here, $|q-1|_3=3$, and, using Lemma \ref{lemma:largeprime},
$q\geq 43$. 
Observe that, since $3$ divides $q-1,$ there is a group of order $3$ normal in a group isomorphic to $\hatt(q-1)^2$ and so, by Lemma
\ref{lemma:normalizersub}, $3\leq \GL$. Thus $b\big|3v$. Now $\Galph$
has the same size as $\Galph$ in Case 1 with $u=6$. The arguments
given there exclude both $b=3v$ and $b=\frac{3}{2}v$ and we are
done.\label{psl353b, psl3FIFTY-FOUR}

Thus we have ruled out all possible actions of line-transitive actions
of $PSL(3,q)$ where $3$ is the unique significant prime.

\section{Case IV: $2|q-1$ is uniquely significant}\label{section:twosig}

In this case $\Galph$ either lies in a parabolic subgroup or in a
subgroup of $PSL(3,q)$ of type 2 or 4. Since $D$ normalizes a Sylow
$2$-subgroup of $PSL(3,q)$, we know that $\Galph$ contains
$D$ for some $\alpha$. Furthermore, by Lemma \ref{lemma:normalizersub}, either
$\Galph\geq \hatt GL(2,q)$ or $D<\GL$.

\subsection{Case 1: $\Galph$ lies inside a group of type 4 only}

In this case $\Galph = PSL(3,q_0)$ or $PSL(3,q_0).3$ for some $q_0$
where $q=q_0^a, a$ odd. Then $D<\GL$ and so $\frac{q-1}{\mu}$ divides into
$3|PSL(3,q_0)|$. We must have $q=q_0^3$. But then $PSL(3,q_0)$ does
not contain an element of order $\frac{q_0^3-1}{\mu}$ and so
$D\not<PSL(3,q_0)$ and this case is also excluded.\label{psl3FIFTY-SIX}

\subsection{Case 2: $\Galph$ lies inside a group of type 2} 

Here $\Galph$ is non-maximal, $q \equiv 1(4)$ and $\Galph$
contains a cyclic subgroup of order $q-1/\mu$. We have two possibilities: 

\begin{enumerate}
\item	$\Galph=A:2$ where $A\leq\hatt (q-1)^2$ and
$|A|=a\frac{q-1}{\mu}$. Then $A$ is proper normal in $\hatt(q-1)^2$
for $a<q-1$ and proper normal in $\hatt(q-1)^2:S_3$ for $a=q-1$. Thus we may conclude, by Lemma \ref{lemma:normalizersub}, that $\GL=A$. We can conclude that $\GL$ contains at most 3 involutions.

\item	We suppose that $3|(q-1)$ and $\Galph = (\frac{q-1}{3}\times\frac{q-1}{3}):S_3$. In this case, $(\frac{q-1}{3}\times\frac{q-1}{3})$ is normal in $\hatt(q-1)^2$ and hence lies in $\GL$. We can conclude that $\ordGL= 3(\frac{q-1}{3})^2$ and $\GL$ contains at most 9 involutions.

\end{enumerate}

Consider the first case. Since $\GL$ contains at most 3 involutions, we may appeal to Lemma \ref{lemma:invineq} to give,
$$k\leq \frac{r_gv}{n_g}+1 = \frac{3q(q+1)(q-1)}{2a}+1.$$\label{psl3FIFTY-SEVENb}
This implies that,
$$k(k-1)< \frac{9}{4a^2}q^3(q+1)^2(q-1).$$
Now we know that $k(k-1) = \frac{v-1}{2}$. Thus
$$\frac{v-1}{2}=\frac{q^3(q^2+q+1)(q+1)(q-1)-2a}{4a}<\frac{9}{4a^2}q^3(q+1)^2(q-1).$$
Hence $q<\frac{9}{a}$ which is impossible.

We move on to the next possibility: $H = (\frac{q-1}{3}\times\frac{q-1}{3})$ lies inside $\GL$ with index 3. Now $H$ contains 3 involutions, hence $\GL$ must contain at most 9 involutions. Once again we appeal to Lemma \ref{lemma:invineq} to give,
$$k\leq \frac{r_gv}{n_g}+1 = \frac{9q(q+1)}{2}+1.$$
This gives,
$$\frac{v-1}{2}=k(k-1)<\frac{41q^2(q+1)^2}{2}.$$\label{psl3FIFTY-EIGHTb}
Given our value for $v$ we may conclude that,
$$q^3(q^2+q+1)(q+1)-2<41q^2(q+1)^2.$$
This is only true for $q\leq 7$ which is impossible.

\subsection{Case 3: $\Galph$ lies in a parabolic subgroup} 

Now, for $P$ a parabolic subgroup, $|G:P|=q^2+q+1$. By Lemma
\ref{lemma:cong} this means that any significant prime must divide
$\frac{1}{2}q(q+1)$. Since $2$ is uniquely significant, we may
conclude that $q\equiv 3 (4)$ and $b\big|\frac{1}{2}(q+1)v$. We write
$\Galph = A.B$ where $A$ is an elementary abelian $p$-group and
$B\leq\hatt GL(2,q)$.

Suppose $q\equiv 3 (8)$. Then, by Lemma \ref{lemma:cong},
 $b=2v$. Then, by Lemma \ref{lemma:bdoublev}, any prime $m$ dividing
 into $v$ must be equivalent to $1(4)$. Since $p\equiv 3 (4)$ we have
 $q^3$ dividing into $\ordGalph$. Thus $A=[q^2]$ and $B\geq \hatt
 SL(2,q)$. However this means that $A.B$ is normal in the full parabolic
 subgroup. Hence, by Lemma \ref{lemma:normalizersub},
 either $\GL\geq\Galph$ (which is impossible) or $\Galph$ is the full
 parabolic subgroup. This case has already been excluded.

Thus $q\equiv 7(8)$ and  $B$ is a subgroup
of $\hatt GL(2,q)$ of type 3 or 5. Consider the case where $B$ is a
subgroup of $\hatt GL(2,q)$ of type 3. We examine the possible
situations here:\label{psl3FIFTY-NINE}

\begin{enumerate}

\item	Suppose that $B$ is maximal in $\hatt GL(2,q)$, i.e. $|B| =
2(q^2-1)/\mu$. Then $B$ acts by conjugation on the non-trivial
elements of $A$ with orbits divisible by $q+1$. Thus $|A|=q^2$ or
1. Since $2$ is uniquely significant, $A<\GL$. This is the same
situation as in Subsections \ref{subsection:qplusoneA} and
\ref{subsection:qplusoneC}; precisely
the same arguments as in those sections allow us to exclude the
situation here.

\item	Suppose that $B$ is non-maximal in $\hatt GL(2,q)$. Then $B$
contains a cyclic group $C$ which is normal in $\hatt (q^2-1)$, hence
lies in $\GL$. Furthermore $|A|\big|\ordGL$ since $2$ is uniquely
significant. Thus $\ordGL=|A|.|C|$ and $\Galph = 2|A|.|C|$ and so
$b=2v$. However in this case, by Lemma \ref{lemma:bdoublev}, any prime
$m$ dividing into $v$ must be equivalent to $1(4)$. Here though
$p\equiv 3 (4)$ and $p$ divides into $v$. This is a
contradiction.
\end{enumerate}

Now consider the possibility that $B$ is of type 5. Since $q\equiv
7(8)$, we must have $q=p^a$ where $a$ is odd and so
$B=\hatt \langle SL(2,q_0),V\rangle$.

Suppose first that $q=q_0$ and so $B\geq\hatt SL(2,q)$ and either $A$
is trivial or $A=[q^2]$.

If $A$ is trivial then either $B\lhd \hatt GL(2,q)$ or $B =
\hatt GL(2,q)$. The first option implies that $\GL\geq \Galph$ (which
is impossible). The latter option is the same as in Subsection
\ref{subsection:qplusoneB}; precisely the same arguments as in that
section allow us to exclude the situation here.

If on the other hand $A$ is non-trivial then $A=[q^2]$ and so $\Galph$
is either the full parabolic subgroup (this possibility is already
excluded) or $\Galph$ is normal in the full parabolic subgroup and
$\GL\geq \Galph$ (which is impossible). Thus both possibilities are
excluded when $q=q_0$. We assume that $q=q_0^a,$ $a$ is odd, $a\geq 3,
p\equiv7(8)$ and $D<\GL$.

Now observe that $A.\langle V\rangle$ is a split extension by Schur-Zassenhaus. So we
can take $V$ to be in $\Galph$. Furthermore $\Galph$ must
contain a conjugate of $D$. Then, since $q\geq q_0^3$,
$\langle V\rangle\cong\frac{q-1}{\mu}$ is $G$-conjugate to $D$. The $G$-conjugates of $D$
split into two conjugacy  
classes inside the parabolic subgroup with centralizers isomorphic to $\hatt
[q]:(q-1)^2$ and $\hatt GL(2,q).$ If we factor out the unipotent
subgroup of the maximal parabolic then we see that, in $\Galph/A$, ${\langle V\rangle A}$ is
centralized by $SL(2,q_0)$ and so $\langle V\rangle$ must be centralized in the maximal
parabolic by $\hatt GL(2,q).$ This means that $\langle V\rangle$ acts by conjugation
on the non-identity elements of $A$ with orbits of size
$\frac{q-1}{\mu}$. In fact $B$ has orbits of length a multiple of
$\frac{(q_0+1)(q-1)}{\mu}$ on the non-trivial elements of $A$. Thus
$|A|=q^2$ or $|A|=1$.\label{psl3SIXTY}

Now note that, since $b\big| \frac{1}{2}v(q+1)q$, we know that
$\frac{2q_0(q_0-1)(q-1)}{\mu}\big|\ordGL$.  Thus $\GL$ lies inside a
subgroup of $PSL(3,q)$ of type 1 or 4. \label{psl3sixty-one}

If $\GL$ lies in a subgroup of $PSL(3,q)$ of type 9 then $\GL=SO(3,q).$
If $A$ is trivial then $\ordGL>\ordGalph$ which is a contradiction. If
$A$ is non-trivial then $q^2$ divides into $\ordGL$ which is a
contradiction. \label{psl3sixty-TWO}

If $\GL$ lies in a subgroup of $PSL(3,q)$ of type 4 then
$\GL=PSL(3,q_1)$ or $PSL(3,q_1).3$. Since $D<\GL$ we must have
$q\leq q_1^2$. But $q=p^a$ where $a$ is odd which is a
contradiction.

Thus $\GL$ lies inside a parabolic subgroup of $PSL(3,q)$. So
$\GL=A_1.B_1$ where $A_1$ is elementary abelian and $B_1\leq \hatt
GL(2,q)$. Then $\frac{2(q_0-1)(q-1)}{\mu}$ divides into
$|B_1|$ and $B_1$ is of type 4, 5, 6 or 7. 

If $B_1$ is of type 5 then we must have $B_1\geq SL(2,q_0)$. Since
$D<A_1.B_1$ we require that $B_1$ contains a cycle of length
$\frac{q-1}{2\mu}$ and so $B_1\geq \langle SL(2,q_0),\frac{q-1}{2\mu}\rangle.$ If $A$
is trivial then $|B_1|\geq \frac{1}{2}\ordGalph$ which is a
contradiction. If $A=[q^2]$ then $A_1$ must be non-trivial and $B_1$
has orbits on the non-trivial elements of $A_1$ of size a multiple of
$\frac{(q_0+1)(q-1)}{\mu}$. Thus $|A_1|=q^2$ and $\ordGL\geq
\frac{1}{2}\ordGalph$. By Lemma \ref{lemma:bdoublev}, $p\equiv 1(4)$
which is a contradiction.\label{PSL3SIXTY-TWOa}

If $B_1$ is of type 4, 6 or 7 then $q_0$ divides into
$|A_1|$ and $\GL=A_1.B_1$ is a split extension. Furthermore $A$ is
trivial since $q^2q_0$ cannot divide into $\ordGL$. 

In the case of types 6 and 7, $B_1$ must centralize $EA_1$ in
$\GL/A_1$ where $E$ is a conjugate of $D$. Thus $E$ has an orbit on
the non-trivial elements of $A_1$ of size a multiple of
$\frac{(q-1)}{\mu}$. Thus $|A_1|\geq q$. But
$\ordGalph<q_0^3\frac{q-1}{\mu}$ and $\ordGL>q\frac{q-1}{\mu}$ which is
impossible.\label{psl3SIXTY-THREE} 

We are left with the possibility that $B_1$ is of type 4 and take $D$
to be in $\GL$. Suppose
first that $DA_1$ is central in $B_1=\GL/A_1$. Since $q+1$
does not divide into $b$, $|B_1|_2\geq 2|(q-1)^2|_2$. This implies
that $D$ is centralized in the full parabolic by $\hatt GL(2,q)$ and
$D$ has orbits on $A_1$ of size a multiple of $\frac{(q-1)}{\mu}$. If,
on the other hand, $DA_1$ is not central in $B_1=\GL/A_1$ then it is
not normal either and $|B_1|$ is divisible by $2(\frac{q-1}{\mu})^2$
Then $\GL$ has orbits on the non-trivial elements of $A_1$
of size a multiple of $\frac{(q-1)}{\mu}$. Thus in either case $|A_1|\geq
q$. But $\ordGalph<q_0^3\frac{q-1}{\mu}$ and $\ordGL>q\frac{q-1}{\mu}$
which is impossible.

This deals with all the cases where $2$ is a uniquely significant prime. We
conclude that $PSL(3,q)$ has no line-transitive actions in this case.

\vspace{5.0pt}
We have now dealt with all possibilities for line-transitive actions
of $PSL(3,q)$ on finite linear spaces. Our proof of Theorem A is complete.

\bibliographystyle{plain}
\bibliography{paper}

\end{document}